\documentclass[10pt,twoside]{article}
\usepackage{mathrsfs}
\usepackage{amsmath}
\usepackage{amssymb}
\usepackage{fancyhdr}
\usepackage{latexsym}
\usepackage{bbding}
\usepackage{mathrsfs}
\usepackage{wasysym}

\usepackage{multicol,graphics}

\newtheorem{theorem}{Theorem}[section]

\newtheorem{remark}[theorem]{Remark}

\numberwithin{equation}{section}

\allowdisplaybreaks

 \makeatletter
\begin{document}
\title{{Logarithmical Blow-up Criteria for the Nematic Liquid Crystal Flows}
\thanks{The author is partially supported by the National Natural Science Foundation of China (11171357).}
}
\author{{\small   Qiao Liu$^{1}$ \thanks{\text{E-mail address}: liuqao2005@163.com.
}}
 {\small\quad and\quad  Jihong Zhao$^{2}$ \thanks{\text{E-mail
address}: zhaojihong2007@yahoo.com.cn.}}
\\
%EndAName
{\small  $^{1}$Department of Mathematics, Hunan Normal University, Changsha, Hunan, 410081,}\\
{\small People's Republic of China}\\
{\small   $^{2}$ College of Science, Northwest A\&F
University, Yangling,}\\
{\small   Shaanxi 712100, People's Republic of China}\\
}
\date{}
\maketitle

\begin{abstract}
We investigate the blow-up criterion for local in time classical
solution of the nematic liquid crystal flows in dimension two and
three. More precisely, $0<T_{*}<+\infty$ is the maximal time
interval if and only if (i) for $n=3$,
\begin{align*}
\int_{0}^{T_{*}}\frac{\|\omega\|_{\dot{B}^{0}_{\infty,\infty}}+\|\nabla
d\|_{\dot{B}^{0}_{\infty,\infty}}^{2}}{\sqrt{1+\text{ln}(e+\|\omega\|_{\dot{B}^{0}_{\infty,\infty}}
+\|\nabla d\|_{\dot{B}^{0}_{\infty,\infty}})}}\text{d}t=\infty,
\end{align*}
or
\begin{align*}
\int_{0}^{T_{*}}\frac{\|\nabla
u\|_{\dot{B}^{-1}_{\infty,\infty}}^{2}+\|\nabla
d\|_{\dot{B}^{0}_{\infty,\infty}}^{2}}{\sqrt{1+\text{ln}(e+\|\nabla
u\|_{\dot{B}^{-1}_{\infty,\infty}} +\|\nabla
d\|_{\dot{B}^{0}_{\infty,\infty}})}}\text{d}t=\infty;
\end{align*}
and (ii) for $n=2$,
\begin{align*}
\int_{0}^{T_{*}}\frac{\|\nabla
d\|_{\dot{B}^{0}_{\infty,\infty}}^{2}}{\sqrt{1+\text{ln}(e +\|\nabla
d\|_{\dot{B}^{0}_{\infty,\infty}})}}\text{d}t=\infty.
\end{align*}
\medskip


\end{abstract}

\section{Introduction}\label{Int}

\noindent

In this paper, we are interested in  the following Cauchy problem of
the flow of the nematic liquid crystal material in $n$-dimensions
($n= 2$ or $3$):
\begin{align}
   \label{eq1.1}
&{\partial_{t}}u-\nu\Delta u +(u\cdot\nabla)u+\nabla{P}=-\lambda\nabla\cdot(\nabla d \odot\nabla d)\quad\text{ in }\mathbb{R}^{n}\times (0,+\infty),\\
%--------------------(eq1.1)------------------------------------
   \label{eq1.2}
&\partial_{t}d+(u\cdot\nabla)d=\gamma(\Delta d+|\nabla d|^{2}d)\quad\quad
  \quad\quad\quad\quad\quad\quad\text{ in }\mathbb{R}^{n}\times (0,+\infty),\\
%--------------------(eq1.2)------------------------------------
   \label{eq1.3}
&\quad\quad\nabla\cdot u=0\quad\quad\quad\quad\quad\quad\quad\quad
\quad\quad\quad\quad\quad\quad\quad\quad\text{ in }\mathbb{R}^{n}\times (0,+\infty),\\
%--------------------(eq1.3)------------------------------------
   \label{eq1.4}
&(u,d)|_{t=0}=(u_{0},d_{0})\quad\quad\quad\quad\quad\quad\quad\quad\quad\quad\quad\quad\quad\quad\text{
in }\mathbb{R}^{n},
%--------------------(eq1.4)------------------------------------
\end{align}
where $u(x,t):\mathbb{R}^{n}\times (0,+\infty)\rightarrow
\mathbb{R}^{n}$ is the unknown velocity field of the flow,
$P(x,t):\mathbb{R}^{n}\times (0,+\infty)\rightarrow \mathbb{R}$ is
the scalar pressure and $d:\mathbb{R}^{n}\times
(0,+\infty)\rightarrow \mathbb{S}^{2}$, the unit sphere in
$\mathbb{R}^{3}$, is the unknown (averaged) macroscopic/continuum
molecule orientation of the nematic liquid crystal flow, $\nabla
\cdot u=0$ represents the incompressible condition, $u_{0}$ is a
given initial velocity with $\nabla \cdot u_{0}=0$ in distribution
sense, $d_{0}: \mathbb{R}^{n}\rightarrow \mathbb{S}^{2}$ is a given
initial liquid crystal orientation field, and $\nu,\lambda,\gamma$
are positive constants. The notation $\nabla d\odot\nabla d$ denotes
the $n\times n$ matrix whose $(i,j)$-th entry is given by
$\partial_{i}d\cdot
\partial_{j}d$ ($1\leq i,j\leq n$). Since the concrete values of the constants $\nu$,
$\lambda$ and $\gamma$ do not play a special role in our discussion,
for simplicity, we assume that they all equal to one throughout this
paper.

The system \eqref{eq1.1}--\eqref{eq1.4}  is a simplified version of
the Ericksen-Leslie model \cite{ER,LE}, which can be viewed as the
incompressible Navier--Stokes equations (the case  $d\equiv 1$, see
\cite{YG,KT,PG}) coupling the heat flow of a harmonic map (the case
$u\equiv 0$, see \cite{CS,PG,W}). Mathematical analysis of the
system \eqref{eq1.1}--\eqref{eq1.4} was initially studied by a
series of papers by Lin \cite{L} and Lin and Liu \cite{LL1,LL2}.
Later on, there are many extensive studies devote to the nematic
liquid crystal flows, see
\cite{JHW,HW,HW1,XLW,LLW,LW,LD,LNW,SL,W,WD} and references therein.
For instance, when the dimension $n=2$, Lin, Lin and Wang \cite{LLW}
established global existence of Leray-Hopf type weak solutions to
\eqref{eq1.1}--\eqref{eq1.4} on bounded domain in $\mathbb{R}^{2}$
under suitable initial and boundary value conditions.  Li and Wang
in \cite{XLW} established the existence of local strong solution
with large initial value and the global strong solution with small
initial value for the initial-boundary value problem of system
\eqref{eq1.1}--\eqref{eq1.4}. Wang in \cite{W} proved that if the
initial data $(u_{0},d_{0})\in BMO^{-1}\times BMO$ is sufficiently
small, then system \eqref{eq1.1}--\eqref{eq1.4} exists a global mild
solution. Lin and Wang \cite{LW} established that when the initial
data $(u_{0},d_{0})$ satisfying $u_{0},\nabla d_{0}\in
L^{n}(\mathbb{R}^{n})$, the solution $(u,d)\in
C([0,T];L^{n}(\mathbb{R}^{n}))\times C([0,T];
\dot{W}^{1,n}(\mathbb{R}^{n},\mathbb{S}^{2}))$ to system
\eqref{eq1.1}--\eqref{eq1.4} is unique.

In the present paper, we are interesting in the short time classical
solution to the system \eqref{eq1.1}--\eqref{eq1.4}. Since the
strong solutions of  the heat flow of  harmonic maps must be blowing
up at finite time \cite{CDY}, we cannot expect that
\eqref{eq1.1}--\eqref{eq1.4} has a global smooth solution with
general initial data. It is well-known that if the initial velocity
$u_{0}\in H^{s}(\mathbb{R}^{n},\mathbb{R}^{n})$ with $\nabla\cdot
u_{0}$ and $d_{0}\in H^{s+1}(\mathbb{R}^{n},\mathbb{S}^{2})$ for
$s\geq n$, then there exists $0<T_{*}<+\infty$ depending only on the
initial value such that the system \eqref{eq1.1}--\eqref{eq1.4} has
a unique local classical solution $(u,d)\in \mathbb{R}^{n}\times
[0,T_{*})$ satisfying (see for example \cite{WD})
\begin{align}\label{eq1.5}
&u\in C([0,T];H^{s}(\mathbb{R}^{n},\mathbb{R}^{n}))\cap
C^{1}([0,T];H^{s-1}(\mathbb{R}^{n},\mathbb{R}^{n}))\quad \text{and
}\nonumber\\
& d\in C([0,T];H^{s+1}(\mathbb{R}^{n},\mathbb{S}^{2}))\cap
C^{1}([0,T];H^{s}(\mathbb{R}^{n},\mathbb{S}^{2}))
\end{align}
%--------------------(eq1.5)------------------------------------
for all $0<T<T_{*}$. Here, we emphasize that  such an existence
theorem gives no indication as to whether solutions actually lose
their regularity or the manner in which they may do so.  Assume that
such $T_{*}$ is the maximum value for \eqref{eq1.5} holds, the
purpose of this paper is to characterize such a $T_{*}$.

For the well-known Navier-Stokes equations with dimension $n\geq 3$,
the Serrin conditions (see \cite{JS,PG}) state that if
$0<T_{*}<\infty$ is the first finite singular time of the smooth
solutions $u$, then $u$ does not belong to the class
$L^{\alpha}(0,T_{*};L^{\beta}(\mathbb{R}^{n}))$ for all
$\frac{2}{\alpha}+\frac{n}{\beta}\leq 1$, $2<\alpha<\infty$,
$n<\beta<\infty$. Beale, Kato and Majida in \cite{BKM} proved that
the vorticity $\omega=\nabla\times u$ does not belong to
$L^{1}(0,T_{*};L^{\infty}(\mathbb{R}^{n}))$ if $T_{*}$ is the first
finite singular time. Later on, Kozono and Taniuchi \cite{KT},
Kozono, Ogawa and Taniuchi \cite{KOT} and Guo and Gala \cite{GG11}
improved the results of \cite{BKM} into BMO and Besov space, more
precisely, if $T_{*}$ is the first singular time, then there hold
\begin{align*}
&\int_{0}^{T_{*}}\|\omega \|_{BMO}\text{d}t=\infty;\\
&\int_{0}^{T_{*}}\frac{\|\omega\|_{\dot{B}^{0}_{\infty,\infty}}}
{\sqrt{1+\text{ln}(1+\|\omega\|_{\dot{B}^{0}_{\infty,\infty}})}}
\text{d}x=\infty,
\end{align*}
where $\dot{B}^{0}_{\infty,\infty}$ denotes the homogeneous Besov
space. On the other hand, as for the heat flow of harmonic maps into
$\mathbb{S}^{2}$, Wang \cite{W} established that for $n\geq 2$, the
condition $\nabla d\in L^{\infty}(0,T;L^{n}(\mathbb{R}^{n}))$
implies that the solution $d$ is regular on $(0,T]$, i.e., $d\in
C^{\infty}((0,T]\times \mathbb{R}^{n})$. For the system
\eqref{eq1.1}--\eqref{eq1.4}, when dimension $n=2$, Lin, Lin and
Wang obtained that the local smooth solution $(u,d)$ to
\eqref{eq1.1}--\eqref{eq1.4} can be continued past any time $T>0$
provided that there holds
\begin{align*}
\int_{0}^{T}\|\nabla d (\cdot, t)\|_{L^{4}}^{4}\text{d}t<\infty.
\end{align*}
Huang and Wang \cite{HW1} established that
\begin{align*}
&\int_{0}^{T_{*}} (\|\omega\|_{L^{\infty}}+\|\nabla
d\|_{L^{\infty}}^{2})\text{d}x=\infty \quad \text{ when dimension }
n=3;\\
&\quad\quad\int_{0}^{T_{*}} \|\nabla
d\|_{L^{\infty}}^{2}\text{d}x=\infty \quad \text{ when dimension }
n=2,
\end{align*}
where $0<T_{*}<\infty$ is the first finite singular time. Motivated
by the above cited papers, the purpose of this paper is to establish
blow-up criteria for local smooth solutions of system
\eqref{eq1.1}--\eqref{eq1.4} in term of the homogeneous Besov
spaces.

Our main results are as follows:

\begin{theorem}\label{thm1.1}
For $n=3$, $u_{0}\in H^{3}(\mathbb{R}^{3},\mathbb{R}^{3})$ with
$\nabla\cdot u_{0}=0$ and $d_{0}\in
H^{4}(\mathbb{R}^{3},\mathbb{S}^{2})$, let $T_{*}>0$ be the maximum
value such that the nematic liquid crystal flow
\eqref{eq1.1}--\eqref{eq1.4} has a unique solution $(u,d)$
satisfying \eqref{eq1.5}. If $T_{*}<+\infty$, then
\begin{align}\label{eq1.6}
\int_{0}^{T_{*}}\frac{\|\omega\|_{\dot{B}^{0}_{\infty,\infty}}+\|\nabla
d\|_{\dot{B}^{0}_{\infty,\infty}}^{2}}{\sqrt{1+\text{ln}(e+\|\omega\|_{\dot{B}^{0}_{\infty,\infty}}
+\|\nabla d\|_{\dot{B}^{0}_{\infty,\infty}})}}\text{d}t=+\infty,
\end{align}
%--------------------(eq1.6)------------------------------------
and
\begin{align}\label{eq1.7}
\int_{0}^{T_{*}}\frac{\|\nabla
u\|_{\dot{B}^{-1}_{\infty,\infty}}^{2}+\|\nabla
d\|_{\dot{B}^{0}_{\infty,\infty}}^{2}}{\sqrt{1+\text{ln}(e+\|\nabla
u\|_{\dot{B}^{-1}_{\infty,\infty}} +\|\nabla
d\|_{\dot{B}^{0}_{\infty,\infty}})}}\text{d}t=+\infty,
\end{align}
%--------------------(eq1.7)------------------------------------
where $\omega:=\nabla\times u$ is the vorticity. In particular, it
holds that
\begin{align*}
\lim\sup_{t\rightarrow
T_{*}}(\|\omega\|_{\dot{B}^{0}_{\infty,\infty}}+\|\nabla
d\|_{\dot{B}^{0}_{\infty,\infty}}^{2})=+\infty;\\
\lim\sup_{t\rightarrow T_{*}}(\|\nabla
u\|_{\dot{B}^{-1}_{\infty,\infty}}+\|\nabla
d\|_{\dot{B}^{0}_{\infty,\infty}})=+\infty.
\end{align*}
\end{theorem}
%--------------------(thm1.1)------------------------------------

\begin{remark}\label{rem1.2}
1.\ By the Sobolev imbedding $L^{\infty}(\mathbb{R}^{3})\subset
\dot{B}^{0}_{\infty,\infty}(\mathbb{R}^{3})$,  it is easy to see
that the condition \eqref{eq1.6} is an extension that of \cite{HW1}.

2.\ Notice that $\nabla u\in
\dot{B}^{-1}_{\infty,\infty}(\mathbb{R}^{3})$ is equivalent to $u\in
\dot{B}^{0}_{\infty,\infty}$, it follows that the condition
\eqref{eq1.7} can be replaced by the following condition:
\begin{align*}
\int_{0}^{T_{*}}\frac{\|
u\|_{\dot{B}^{0}_{\infty,\infty}}^{2}+\|\nabla
d\|_{\dot{B}^{0}_{\infty,\infty}}^{2}}{\sqrt{1+\text{ln}(e+\|
u\|_{\dot{B}^{0}_{\infty,\infty}} +\|\nabla
d\|_{\dot{B}^{0}_{\infty,\infty}})}}\text{d}t=+\infty.
\end{align*}
In particular, there holds
\begin{align*}
\lim\sup_{t\rightarrow T_{*}}(\|
u\|_{\dot{B}^{0}_{\infty,\infty}}+\|\nabla
d\|_{\dot{B}^{0}_{\infty,\infty}})=+\infty.
\end{align*}
\end{remark}
%--------------------(rem1.2)------------------------------------

As a byproduct of our proof of Theorem 1.1, we obtain the following
corresponding criterion in dimension two. More precisely, we have

\begin{theorem}\label{cor1.3}
For n=2, $u_{0}\in H^{2}(\mathbb{R}^{2},\mathbb{R}^{2})$ with
$\nabla\cdot u_{0}=0$ and $d_{0}\in
H^{3}(\mathbb{R}^{2},\mathbb{S}^{2})$, let $T_{*}>0$ be the maximum
value such that the nematic liquid crystal flow
\eqref{eq1.1}--\eqref{eq1.4} has a unique solution $(u,d)$
satisfying \eqref{eq1.5}. If $T_{*}<+\infty$, then
\begin{align}\label{eq1.8}
\int_{0}^{T_{*}}\frac{\|\nabla
d\|_{\dot{B}^{0}_{\infty,\infty}}^{2}}{\sqrt{1+\text{ln}(e +\|\nabla
d\|_{\dot{B}^{0}_{\infty,\infty}})}}\text{d}t=+\infty.
\end{align}
%--------------------(eq1.8)------------------------------------
In particular, there holds
\begin{align*}
\lim\sup_{t\rightarrow T_{*}}\|\nabla
d\|_{\dot{B}^{0}_{\infty,\infty}}=+\infty.
\end{align*}
\end{theorem}
%--------------------(cor1.3)------------------------------------

The remaining of the paper is written as follows. Section 2 is
devoted to the proof of Theorem \ref{thm1.1}. In Section 3, we prove
Theorem \ref{cor1.3}. Throughout the paper, C denotes a constant and
may change from line to line; $\|\cdot\|_{X}$ denotes the norm of
space $X(\mathbb{R}^{3})$ or $X(\mathbb{R}^{2})$.

\section{The proof of Theorem \ref{thm1.1}}

In this Section, we shall give the proof of Theorem \ref{thm1.1}.
Before going to the proof, we first review the following two
inequalities, the first one can be found in \cite{KOT} and the
second one can be found in \cite{GG11}:
\begin{align}\label{eq2.1}
\|f\|_{L^{\infty}}\leq
C(1+\|f\|_{\dot{B}^{0}_{\infty,\infty}}\text{ln}^{\frac{1}{2}}(1+\|f\|_{H^{s-1}}))
\end{align}
%--------------------(eq2.1)------------------------------------
for $f\in H^{s-1}(\mathbb{R}^{n})$ with $s>\frac{n}{2}+1$ and $n\geq
2$.
\begin{align}\label{eq2.2}
\|f\|_{L^{4}}\leq
\|f\|_{\dot{B}^{-1}_{\infty,\infty}}^{\frac{1}{2}}\|\nabla
f\|_{L^{2}}^{\frac{1}{2}}
\end{align}
%--------------------(eq2.2)------------------------------------
for $f\in \dot{H}^{1}(\mathbb{R}^{3})\cap
\dot{B}^{-1}_{\infty,\infty}(\mathbb{R}^{3})$.

\textit{Case I:} We now give the proof of \eqref{eq1.6} under the
assumptions of Theorem 1.1.

Since we deal with the local smooth solutions, and notice that
$[0,T_{*})$ is the maximal existence interval of local smooth
solution associated with initial value ($u_{0},d_{0}$). We prove
Theorem 1.1 arguing by contradiction. Suppose, that \eqref{eq1.6} is
not true. Then there is $0<M<\infty$ such that
\begin{align}\label{eq2.3}
\int_{0}^{T_{*}}\frac{\|\omega\|_{\dot{B}^{0}_{\infty,\infty}}+\|\nabla
d\|_{\dot{B}^{0}_{\infty,\infty}}^{2}}{\sqrt{1+\text{ln}(e+\|\omega\|_{\dot{B}^{0}_{\infty,\infty}}
+\|\nabla d\|_{\dot{B}^{0}_{\infty,\infty}})}}\text{d}t\leq M.
\end{align}
%--------------------(eq2.3)------------------------------------
We will show that if assumption \eqref{eq2.3} holds, then there
holds
\begin{align}\label{eq2.4}
\lim_{t\rightarrow T_{*}^{-}}(\|\nabla^{3}
u(t,\cdot)\|_{L^2}^{2}+\|\nabla^{4} d(t,\cdot)\|_{L^2}^{2})\leq C,
\end{align}
%--------------------(eq2.4)------------------------------------
for some positive constant $C$ depends only on $u_{0}, d_{0}, T_{*}$
and $M$. The estimate \eqref{eq2.4} is enough to extend the smooth
solution ($u$, $d$) beyond to $ T_{*}$. That is to say, $[0,T_{*})$
is not a maximal interval of existence, which leads to the
contradiction.

We first taking $\nabla \times$ on \eqref{eq1.1}, it follows that
\begin{align}\label{eq2.5}
\omega_{t}-\Delta \omega +u\cdot\nabla \omega=\omega\cdot\nabla
u-\nabla\times(\Delta d\cdot \nabla d),
\end{align}
%--------------------(eq2.5)------------------------------------
where we have used the facts that $\nabla\cdot(\nabla d\otimes
\nabla d)=\nabla(\frac{|\nabla d|^{2}}{2})+\Delta d\cdot\nabla d$
and $\nabla\times \nabla(\frac{|\nabla d|^{2}}{2})=0$. Multiplying
\eqref{eq2.5} with $\omega$ and integrating over $\mathbb{R}^{3}$,
we obtain
\begin{align}\label{eq2.6}
\frac{1}{2}\frac{d}{dt}\|\omega (\cdot,t)\|_{L^{2}}^{2}+\|\nabla
\omega\|_{L^{2}}^{2}=\int_{\mathbb{R}^{3}}[(\omega\cdot\nabla)u\cdot
\omega+(\Delta d\cdot\nabla d)\cdot\nabla\times
\omega]\text{d}x:=I_{1}+I_{2},
\end{align}
%--------------------(eq2.6)------------------------------------
where we have used the fact that $\operatorname{div} u=0$ implies
that $\int_{\mathbb{R}^{3}}(u\cdot\nabla)\omega\cdot\omega\text{d}x
=\frac{1}{2}\int_{\mathbb{R}^{3}}(u\cdot\nabla)|\omega|^{2}\text{d}x=0$.
By using the H\"{o}lder's inequality and \eqref{eq2.1} with $s=3$,
we can estimate $I_{1}$ as
\begin{align}\label{eq2.7}
  I_{1}\leq& C\|\omega \|_{L^{\infty}}\|\nabla
  u\|_{L^{2}}\|\omega\|_{L^{2}}\leq
  C\|\omega\|_{L^{\infty}}\|\omega\|_{L^{2}}^{2}\nonumber\\
  \leq&
  C(1+\|\omega\|_{\dot{B}^{0}_{\infty,\infty}}\text{ln}^{\frac{1}{2}}(e+\|\Lambda^{2}
  \omega\|_{L^{2}}))\|\omega\|_{L^{2}}^{2}\nonumber\\
  \leq&
  C(1+\|\omega\|_{\dot{B}^{0}_{\infty,\infty}}\text{ln}^{\frac{1}{2}}(e+\|\Lambda^{3}
  u\|_{L^{2}}))\|\omega\|_{L^{2}}^{2}\nonumber\\
  \leq&
  C\frac{\|\omega\|_{\dot{B}^{0}_{\infty,\infty}}}{\sqrt{1
  +\text{ln}(e+\|\omega\|_{\dot{B}^{0}_{\infty,\infty}}
  +\|\nabla d\|_{\dot{B}^{0}_{\infty,\infty}})}}\text{ln}^{\frac{1}{2}}(e+\|\Lambda^{3}
  u\|_{L^{2}})\nonumber\\
  &\times(1
  +\text{ln}(e+\|\omega\|_{\dot{B}^{0}_{\infty,\infty}}
  +\|\nabla d\|_{\dot{B}^{0}_{\infty,\infty}}))^{\frac{1}{2}}\|\omega\|_{L^{2}}^{2}+C\|\omega\|_{L^{2}}^{2}\nonumber\\
  \leq&
  C\frac{\|\omega\|_{\dot{B}^{0}_{\infty,\infty}}}{\sqrt{1
  +\text{ln}(e+\|\omega\|_{\dot{B}^{0}_{\infty,\infty}}
  +\|\nabla d\|_{\dot{B}^{0}_{\infty,\infty}})}}\text{ln}^{\frac{1}{2}}(e+\|\Lambda^{3}
  u\|_{L^{2}}+\|\Lambda^{4} d\|_{L^{2}})\nonumber\\
  &\times(1
  +\text{ln}(e+\|\Lambda^{2}\omega\|_{L^{2}}
  +\|\Lambda^{3}\nabla d\|_{L^{2}}))^{\frac{1}{2}}\|\omega\|_{L^{2}}^{2}+C\|\omega\|_{L^{2}}^{2}\nonumber\\
  \leq& C\frac{\|\omega\|_{\dot{B}^{0}_{\infty,\infty}}}{\sqrt{1
  +\text{ln}(e+\|\omega\|_{\dot{B}^{0}_{\infty,\infty}}
  +\|\nabla d\|_{\dot{B}^{0}_{\infty,\infty}})}}\text{ln}(e+\|\Lambda^{3}
  u\|_{L^{2}}+\|\Lambda^{4}
  d\|_{L^{2}})\|\omega\|_{L^{2}}^{2}+C\|\omega\|_{L^{2}}^{2},
\end{align}
%--------------------(eq2.7)------------------------------------
where $\Lambda: =(-\Delta)^{\frac{1}{2}}$, and we have used the fact
\begin{align*}
\nabla u=(-\Delta)^{-1}\nabla (\nabla\times\omega) \text{ implies }
\|\nabla u\|_{L^{2}}\leq C\|\omega\|_{L^{2}},
\end{align*}
and the  following standard Sobolev imbedding
\begin{align}\label{eq2.8}
H^{3}(\mathbb{R}^{3})\subseteq H^{2}(\mathbb{R}^{3})\subseteq
L^{\infty}(\mathbb{R}^{3})\subseteq BMO(\mathbb{R}^{3})\subseteq
\dot{B}^{0}_{\infty,\infty}(\mathbb{R}^{3}).
\end{align}
%--------------------(eq2.8)------------------------------------
As for $I_{2}$, by using the H\"{o}lder's inequality and
\eqref{eq2.1} with $s=4$, we get
\begin{align}\label{eq2.9}
I_{2}\leq& C\||\Delta d||\nabla d|\|_{L^{2}}\|\nabla
\omega\|_{L^{2}}\leq C \|\nabla d\|_{L^{\infty}}\|\Delta
d\|_{L^{2}}\|\nabla \omega\|_{L^{2}}\nonumber\\
\leq& \frac{1}{4}\|\nabla \omega\|_{L^{2}}^{2}+C\|\nabla
d\|_{L^{\infty}}^{2}\|\Delta d\|_{L^{2}}^{2}\nonumber\\
\leq& \frac{1}{4}\|\nabla \omega\|_{L^{2}}^{2}+C(1+\|\nabla
d\|_{\dot{B}^{0}_{\infty,\infty}}^{2}\text{ln}(e+\|\Lambda^{3}\nabla
d \|_{L^{2}}))\|\Delta d\|_{L^{2}}^{2}\nonumber\\
\leq& \frac{1}{4}\|\nabla \omega\|_{L^{2}}^{2}+C\frac{\|\nabla
d\|_{\dot{B}^{0}_{\infty,\infty}}^{2}}{\sqrt{1+\text{ln}
(e+\|\omega\|_{\dot{B}^{0}_{\infty,\infty}}+\|\nabla
d\|_{\dot{B}^{0}_{\infty,\infty}})}}\text{ln}(e+\|\Lambda^{3}\nabla
d \|_{L^{2}})\nonumber\\
&\times (1+\text{ln}
(e+\|\omega\|_{\dot{B}^{0}_{\infty,\infty}}+\|\nabla
d\|_{\dot{B}^{0}_{\infty,\infty}}))^{\frac{1}{2}}\|\Delta
d\|_{L^{2}}^{2}+C\|\Delta d\|_{L^{2}}^{2}\nonumber\\
\leq& \frac{1}{4}\|\nabla \omega\|_{L^{2}}^{2}+C\frac{\|\nabla
d\|_{\dot{B}^{0}_{\infty,\infty}}^{2}}{\sqrt{1+\text{ln}
(e+\|\omega\|_{\dot{B}^{0}_{\infty,\infty}}+\|\nabla
d\|_{\dot{B}^{0}_{\infty,\infty}})}}[\text{ln}(e+\|\Lambda^{3}u\|_{L^{2}}+\|\Lambda^{4}
d \|_{L^{2}})]^{\frac{3}{2}}\|\Delta d\|_{L^{2}}^{2}\nonumber\\
&+C\|\Delta d\|_{L^{2}}^{2}.
\end{align}
%--------------------(eq2.9)------------------------------------
Inserting estimates \eqref{eq2.7} and \eqref{eq2.9} into
\eqref{eq2.6}, it follows that
\begin{align}\label{eq2.10}
&\frac{d}{dt}\|\omega (\cdot,t)\|_{L^{2}}^{2}+\frac{3}{2}\|\nabla
\omega\|_{L^{2}}^{2}\leq C(\|\omega\|_{L^{2}}^{2}+\|\Delta
d\|_{L^{2}}^{2})\nonumber\\
&+C\frac{\|\nabla
d\|_{\dot{B}^{0}_{\infty,\infty}}^{2}}{\sqrt{1+\text{ln}
(e+\|\omega\|_{\dot{B}^{0}_{\infty,\infty}}+\|\nabla
d\|_{\dot{B}^{0}_{\infty,\infty}})}}[\text{ln}(e+\|\Lambda^{3}u\|_{L^{2}}+\|\Lambda^{4}
d \|_{L^{2}})]^{\frac{3}{2}}(\|\omega\|_{L^{2}}^{2}+\|\Delta
d\|_{L^{2}}^{2}).
\end{align}
%--------------------(eq2.10)------------------------------------

Taking $\Delta $ on equation \eqref{eq1.3}, multiplying $\Delta d$
and integrating over $\mathbb{R}^{3}$, one obtains
\begin{align}\label{eq2.11}
\frac{1}{2}\frac{d}{dt}\|\Delta d (\cdot,t)\|_{L^{2}}^{2}+\|\nabla
\Delta d\|_{L^{2}}^{2}=-\!\int_{\mathbb{R}^{3}}\Delta(u\cdot\nabla
d)\cdot\Delta d\text{d}x+ \!\int_{\mathbb{R}^{3}}\Delta(|\nabla
d|^{2}d)\cdot\Delta d\text{d}x:=I_{3}+I_{4}.
\end{align}
%--------------------(eq2.11)------------------------------------
Notice that the fact $\operatorname{div} u=0$ implies that
\begin{align*}
\int_{\mathbb{R}^{3}}(u\cdot\nabla)\Delta d\cdot\Delta
d\text{d}x=\frac{1}{2}\int_{\mathbb{R}^{3}}u\cdot\nabla |\Delta
d|^{2}\text{d}x=0,
\end{align*}
and the equality
\begin{align*}
\nabla\times \omega=\nabla\times(\nabla\times u)=\nabla(\nabla\cdot
u)-\Delta u=-\Delta u\quad\text{implies } \|\Delta u\|_{L^{2}}\leq
C\|\nabla \omega\|_{L^{2}}.
\end{align*}
Hence we can estimate $I_{3}$ as
\begin{align}\label{eq2.12}
I_{3}\leq& \int_{\mathbb{R}^{3}}|\Delta u||\nabla d||\Delta
d|\text{d}x+2\int_{\mathbb{R}^{3}}|\nabla u||\nabla^{2}d||\Delta
d|\text{d}x\nonumber\\
\leq& \|\Delta u\|_{L^{2}}\|\nabla d\|_{L^{\infty}}\|\Delta
d\|_{L^{2}}+\|\nabla u\|_{L^{2}}\|\nabla^{2}d\|_{L^{4}}\|\Delta
d\|_{L^{4}}\nonumber\\
\leq& \|\nabla \omega\|_{L^{2}}\|\nabla d\|_{L^{\infty}}\|\Delta
d\|_{L^{2}}+\|\omega\|_{L^{2}}\|\nabla
d\|_{L^{\infty}}\|\nabla\Delta
d\|_{L^{2}}\nonumber\\
\leq& \frac{1}{4}(\|\nabla \omega\|_{L^{2}}^{2}+\|\nabla \Delta
d\|_{L^{2}}^{2})+C\|\nabla
d\|_{L^{\infty}}^{2}(\|\omega\|_{L^{2}}^{2}+\|\Delta
d\|_{L^{2}}^{2})\nonumber\\
\leq& \frac{1}{4}(\|\nabla \omega\|_{L^{2}}^{2}+\|\nabla \Delta
d\|_{L^{2}}^{2})+C(1+\|\nabla
d\|_{\dot{B}^{0}_{\infty,\infty}}^{2}\text{ln}(e+\|\Lambda^{3}\nabla
d\|_{L^{2}}))(\|\omega\|_{L^{2}}^{2}+\|\Delta
d\|_{L^{2}}^{2})\nonumber\\
\leq& \frac{1}{4}(\|\nabla \omega\|_{L^{2}}^{2}+\|\nabla \Delta
d\|_{L^{2}}^{2})+C(\|\omega\|_{L^{2}}^{2}+\|\Delta
d\|_{L^{2}}^{2})+\nonumber\\
&\! C\!\!\frac{\|\nabla
d\|_{\dot{B}^{0}_{\infty,\infty}}^{2}}{\sqrt{\!1\!+\text{ln}
(e\!+\|\omega\|_{\dot{B}^{0}_{\infty,\infty}}\!\!+\|\nabla
d\|_{\dot{B}^{0}_{\infty,\infty}})}}[\text{ln}(e\!+\|\Lambda^{3}u\|_{L^{2}}\!+\|\Lambda^{4}
d \|_{L^{2}})]^{\frac{3}{2}}(\|\omega\|_{L^{2}}^{2}\!+\|\Delta
d\|_{L^{2}}^{2}),
\end{align}
%--------------------(eq2.12)------------------------------------
where we have used the Gagliardo-Nirenberg inequalities
\begin{align*}
\|\nabla^{2}d\|_{L^{4}}\leq \|\nabla
d\|_{L^{\infty}}^{\frac{1}{2}}\|\nabla \Delta
d\|_{L^{2}}^{\frac{1}{2}}\quad\text{and}\quad\|\Delta
d\|_{L^{4}}\leq \|\nabla d\|_{L^{\infty}}^{\frac{1}{2}}\|\nabla
\Delta d\|_{L^{2}}^{\frac{1}{2}}.
\end{align*}
%for all $f\in \dot{W}^{1,\infty}(\mathbb{R}^{3})\cap
%\dot{H}^{3}(\mathbb{R}^{3})$.
To estimate $I_{4}$, notice that the
condition $|d|=1$,
\begin{align*}
I_{4}=&\int_{\mathbb{R}^{3}}[\Delta(|\nabla d|^{2})d\cdot\Delta
d+2\nabla(|\nabla d|^{2})\nabla d\cdot\Delta d+|\nabla d|^{2}|\Delta
d|^{2}]\text{d}x\nonumber\\
= &\int_{\mathbb{R}^{3}}[-\nabla(|\nabla d|^{2})\nabla d\cdot\Delta
d-\nabla(|\nabla d|^{2}) d\cdot\nabla\Delta d+2\nabla(|\nabla
d|^{2})\nabla d\cdot\Delta d+|\nabla d|^{2}|\Delta
d|^{2}]\text{d}x\nonumber\\
\leq& \int_{\mathbb{R}^{3}}[6|\nabla d|^{2}|\nabla^{2} d||\Delta
d|+2|\nabla d||\nabla^{2} d||\nabla\Delta d|+|\nabla d|^{2}|\Delta
d|^{2}]\text{d}x\nonumber\\
\leq & C( \|\nabla d\|_{L^{\infty}}^{2}\|\Delta
d\|_{L^{2}}\|\nabla^{2} d\|_{L^{2}}+\|\nabla
d\|_{L^{\infty}}\|\nabla^{2}d\|_{L^{2}}\|\nabla\Delta
d\|_{L^{2}}+\|\nabla d\|_{L^{\infty}}^{2}\|\Delta d\|_{L^{2}}^{2}).
\end{align*}
By the standard calculation, it is easy to obtain the equality
$\|\nabla^{2} d\|_{L^{2}}=\|\Delta d\|_{L^{2}}$. Hence,
\begin{align}\label{eq2.13}
I_{4}\leq& C(\|\nabla d\|_{L^{\infty}}^{2}\|\Delta d\|_{L^{2}}^{2}
+\|\nabla d\|_{L^{\infty}}\|\nabla^{2}d\|_{L^{2}}\|\nabla\Delta
d\|_{L^{2}})\nonumber\\
\leq& \frac{1}{4}\|\nabla \Delta d\|_{L^{2}}^{2}+C\|\nabla
d\|_{L^{\infty}}^{2}\|\Delta d\|_{L^{2}}^{2}\nonumber\\
\leq& \frac{1}{4}\|\nabla \Delta d\|_{L^{2}}^{2}+C(1+\|\nabla
d\|_{\dot{B}^{0}_{\infty,\infty}}^{2}\text{ln}(e+\|\Lambda^{3}\nabla
d\|_{L^{2}}))\|\Delta
d\|_{L^{2}}^{2} \nonumber\\
\leq& \frac{1}{4}\|\nabla \Delta d\|_{L^{2}}^{2}+C\|\Delta
d\|_{L^{2}}^{2}+\nonumber\\
&\!C\!\!\frac{\|\nabla
d\|_{\dot{B}^{0}_{\infty,\infty}}^{2}}{\sqrt{\!1\!+\text{ln}
(e\!+\|\omega\|_{\dot{B}^{0}_{\infty,\infty}}\!\!+\|\nabla
d\|_{\dot{B}^{0}_{\infty,\infty}})}}[\text{ln}(e\!+\|\Lambda^{3}u\|_{L^{2}}\!+\|\Lambda^{4}
d \|_{L^{2}})]^{\frac{3}{2}}(\|\omega\|_{L^{2}}^{2}\!+\|\Delta
d\|_{L^{2}}^{2}).
\end{align}
%--------------------(eq2.13)------------------------------------
Inserting estimates \eqref{eq2.12} and \eqref{eq2.13} into
\eqref{eq2.11}, it follows that
\begin{align}\label{eq2.14}
&\frac{d}{dt}\|\Delta d (\cdot,t)\|_{L^{2}}^{2}+\|\nabla \Delta
d\|_{L^{2}}^{2}\leq \frac{1}{2}\|\nabla
\omega\|_{L^{2}}^{2}+C(\|\omega\|_{L^{2}}^{2}+\|\Delta
d\|_{L^{2}}^{2})\nonumber\\
&\! +C\!\!\frac{\|\nabla
d\|_{\dot{B}^{0}_{\infty,\infty}}^{2}}{\sqrt{\!1\!+\text{ln}
(e\!+\|\omega\|_{\dot{B}^{0}_{\infty,\infty}}\!\!+\|\nabla
d\|_{\dot{B}^{0}_{\infty,\infty}})}}[\text{ln}(e\!+\|\Lambda^{3}u\|_{L^{2}}\!+\|\Lambda^{4}
d \|_{L^{2}})]^{\frac{3}{2}}(\|\omega\|_{L^{2}}^{2}\!+\|\Delta
d\|_{L^{2}}^{2}).
\end{align}
%--------------------(eq2.14)------------------------------------

Due to \eqref{eq2.3}, one concludes that for any small constant
$\varepsilon>0$, there exists $T=T(\varepsilon)\in(0,T_{*})$ such
that
\begin{align*}
\int_{T}^{T_{*}}\frac{\|\omega\|_{\dot{B}^{0}_{\infty,\infty}}+\|\nabla
d\|_{\dot{B}^{0}_{\infty,\infty}}^{2}}{\sqrt{1+\text{ln}(e+\|\omega\|_{\dot{B}^{0}_{\infty,\infty}}
+\|\nabla d\|_{\dot{B}^{0}_{\infty,\infty}})}}\text{d}t\leq
\varepsilon.
\end{align*}
For any $T<t\leq T_{*}$, we let
\begin{align}\label{eq2.15}
y(t):=\sup_{T<\tau\leq t} (\|\Lambda^{3}
u\|_{L^{2}}^{2}+\|\Lambda^{4} d\|_{L^{2}}^{2}).
\end{align}
%--------------------(eq2.15)------------------------------------
Then, by abbreviately denoting $C\{\|\omega
(\cdot,T(\varepsilon))\|_{L^{2}}^{2}+\|\Delta
d(\cdot,T(\varepsilon))\|_{L^{2}}^{2}\}$ as $C(\varepsilon)$, and
putting \eqref{eq2.10}, \eqref{eq2.14} and \eqref{eq2.15} together,
we find that
\begin{align}\label{eq2.16}
\frac{d}{dt}&(\|\omega (\cdot,t)\|_{L^{2}}^{2}+\|\Delta d
(\cdot,t)\|_{L^{2}}^{2})+\|\nabla \omega\|_{L^{2}}^{2}+\|\nabla
\Delta d\|_{L^{2}}^{2} \leq C(\|\omega\|_{L^{2}}^{2}+\|\Delta
d\|_{L^{2}}^{2})\nonumber\\
&+C\frac{\|\nabla
d\|_{\dot{B}^{0}_{\infty,\infty}}^{2}}{\sqrt{\!1\!+\text{ln}
(e\!+\|\omega\|_{\dot{B}^{0}_{\infty,\infty}}\!\!+\|\nabla
d\|_{\dot{B}^{0}_{\infty,\infty}})}}
[\text{ln}(e+y(t))]^{\frac{3}{2}} (\|\omega\|_{L^{2}}^{2}+\|\Delta
d\|_{L^{2}}^{2}).
\end{align}
%--------------------(eq2.16)------------------------------------
By Gronwall's inequality to \eqref{eq2.16} in the interval  $ [T,
t)$, we have
\begin{align}\label{eq2.17}
&\|\omega(\cdot,t)\|_{L^{2}}^{2}+\|\Delta d (\cdot,
t)\|_{L^{2}}^{2}\nonumber\\
 \leq& (\|\omega
(\!\cdot,\!T)\|_{\!L^{\!2}}^{\!2}\!\!+\|\Delta d
(\!\cdot,\!T)\|_{\!L^{\!2}}^{\!2})\exp\!\!\left\{\!\!\int_{T}^{t}\!C\text{d}\tau\!+C\!\!\int_{T}^{t}\!\!
\frac{\|\omega\|_{\dot{B}^{0}_{\infty,\infty}}+\|\nabla
d\|_{\dot{B}^{0}_{\infty,\infty}}^{2}}{\!\sqrt{1\!+\!\text{ln}(e\!+\!\|\omega\|_{\!\dot{B}^{0}_{\!\infty,\!\infty}}
\!\!+\!\|\nabla d\|_{\!\dot{B}^{0}_{\!\infty,\!\infty}})}}\text{d}\tau\text{ln}^{\frac{3}{2}}(e\!+y(t))\!\!\right\}\nonumber\\
 \leq& (\|\omega
(\cdot,T)\|_{L^{2}}^{2}+\|\Delta d
(\cdot,T)\|_{L^{2}}^{2})\exp\left\{C[(T_{*}-T)+\varepsilon\text{ln}^{\frac{3}{2}}(e+y(t))]\right\}\nonumber\\
\leq &
C(\varepsilon)\exp(C\varepsilon(1+\text{ln}^{\frac{3}{2}}(e+y(t))))
\leq C(\varepsilon) (e+y(t))^{\frac{3C\varepsilon}{2}},
\end{align}
%--------------------(eq2.17)------------------------------------
where $C$ is the positive constant whose value is independent of
either $\epsilon$ or $T$, and $C(\varepsilon)$ is a bounded positive
constant depending on $\epsilon$ which may change from line to line.

Next, we will derive the estimate of  $y(t)$ defined by
\eqref{eq2.15}. To this end, we need to introduce the following
commutator and product estimates (see \cite{KP}):
\begin{align}
   \label{eq2.18}
&\|\Lambda^{\alpha}(fg)-f\Lambda^{\alpha}g\|_{L^{p}}\leq C (\|\nabla
f\|_{L^{p_{1}}}\|\Lambda^{\alpha-1}g\|_{L^{q_{1}}}+\|\Lambda^{\alpha}f\|_{L^{p_{2}}}\|g\|_{L^{q_{2}}});\\
%--------------------(eq2.18)------------------------------------
   \label{eq2.19}
&\|\Lambda^{\alpha}(fg)\|_{L^{p}}\leq C
(\|f\|_{L^{p_{1}}}\|\Lambda^{\alpha}g\|_{L^{q_{2}}}+\|\Lambda^{\alpha}f\|_{L^{p_{2}}}\|g\|_{L^{q_{2}}})
\end{align}
%--------------------(eq2.19)------------------------------------
with $\alpha>0$, $1<p,p_{1},p_{2},q_{1},q_{2}<\infty$ and
$\frac{1}{p}=\frac{1}{p_{1}}+\frac{1}{q_{1}}=\frac{1}{p_{2}}+\frac{1}{q_{2}}$.

Applying $\Lambda^{3}$ on \eqref{eq1.1}, multiplying $\Lambda^{3}u$
and integrating over $\mathbb{R}^{3}$, one obtains
\begin{align}\label{eq2.20}
\!\frac{1}{2}\frac{d}{dt}\|\Lambda^{3}u
(\cdot,t)\|_{\!L^{\!2}}^{2}\!+\!\|\Lambda^{4}u
(\cdot,t)\|_{\!L^{\!2}}^{2}\!=\!\!-\!\!\int_{\mathbb{R}^{3}}\!\!
\Lambda^{3}(u\cdot\nabla
u)\cdot\Lambda^{3}u\text{d}x\!-\!\!\!\int_{\mathbb{R}^{3}}\!\!\Lambda^{3}(\Delta
d\cdot\nabla
d)\cdot\Lambda^{3}u\text{d}x\!:=\!I_{\!5}\!\!+\!I_{\!6},
\end{align}
%--------------------(eq2.20)------------------------------------
where we have used the fact that $\operatorname{div} u=0$ implies
$\int_{\mathbb{R}^{3}}\Lambda^{3}\nabla(\frac{|\nabla
d|^{2}}{2})\cdot\Lambda^{3} u\text{d}x=0$. Applying \eqref{eq2.18},
it follows that
\begin{align}\label{eq2.21}
I_{5}=&\int_{\mathbb{R}^{3}}[\Lambda^{3}(u\cdot\nabla
u)-u\cdot\nabla\Lambda^{3}u]\cdot\Lambda^{3}u\text{d}x\nonumber\\
\leq& C\|[\Lambda^{3}(u\cdot\nabla
u)-u\cdot\nabla\Lambda^{3}u]\|_{L^{\frac{3}{2}}}\|\Lambda^{3}u\|_{L^{3}}\nonumber\\
\leq &C\|\nabla u\|_{L^{3}}\|\Lambda^{3} u\|_{L^{3}}^{2}\leq
C\|\nabla
u\|_{L^{2}}^{\frac{13}{12}}\|\Lambda^{3}u\|_{L^{2}}^{\frac{1}{4}}\|\Lambda^{4}u\|_{L^{2}}^{\frac{5}{3}}\nonumber\\
\leq& \frac{1}{4}\|\Lambda^{4} u\|_{L^{2}}^{2}+C\|\nabla
u\|_{L^{2}}^{\frac{13}{2}}\|\Lambda^{3}
u\|_{L^{2}}^{\frac{3}{2}}\nonumber\\
\leq& \frac{1}{4}\|\Lambda^{4}
u\|_{L^{2}}^{2}+C_{0}\|\omega\|_{L^{2}}^{\frac{13}{2}}\|\Lambda^{3}
u\|_{L^{2}}^{\frac{3}{2}}.
\end{align}
%--------------------(eq2.21)------------------------------------
Here we have used the following Gagliardo-Nirenberg inequalities:
\begin{align*}
\|\nabla u\|_{L^{3}}\leq C\|\nabla
u\|_{L^{2}}^{\frac{3}{4}}\|\Lambda^{3} u\|_{L^{2}}^{\frac{1}{4}}
\text{ and }
 \|\Lambda^{3} u\|_{L^{3}}\leq C\|\nabla
u\|_{L^{2}}^{\frac{1}{6}}\|\Lambda^{4}u\|_{L^{2}}^{\frac{5}{6}}.
\end{align*}
For $I_{6}$, applying the H\"{o}lder's inequality and the Leibniz's
rule,  we have
\begin{align}\label{eq2.22}
I_{6}= &\int_{\mathbb{R}^{3}}\Lambda^{2}(\Delta d\cdot\nabla d)\cdot
\Lambda^{4}u\text{d}x\nonumber\\
\leq & \frac{1}{4}\|\nabla^{4}
u\|_{L^{2}}^{2}+C\int_{\mathbb{R}^{3}}|\Lambda^{2}(\Delta
d\cdot\nabla d)|^{2}\text{d}x\nonumber\\
\leq &\frac{1}{4}\|\nabla^{4}
u\|_{L^{2}}^{2}+C\int_{\mathbb{R}^{3}}(|\Lambda^{4}d|^{2}|\nabla
d|^{2}+|\Lambda^{2} d|^{2}|\Lambda^{3} d|^{2})\text{d}x\nonumber\\
\leq& \frac{1}{4}\|\nabla^{4} u\|_{L^{2}}^{2}+C(\|\nabla
d\|_{L^{6}}^{2}\|\Lambda^{4} d\|_{L^{3}}^{2}+\|\Lambda^{2}d\|_{L^{4}}^{2}\|\Lambda^{3}d\|_{L^{4}}^{2})\nonumber\\
\leq& \frac{1}{4}\|\nabla^{4} u\|_{L^{2}}^{2}+C(\|\Delta
d\|_{L^{2}}^{\frac{7}{3}}\|\Lambda^{5}d\|_{L^{2}}^{\frac{5}{3}}+
\|\Delta d\|_{L^{2}}^{\frac{19}{6}}\|\Lambda^{5}
d\|_{L^{2}}^{\frac{5}{6}})\nonumber\\
\leq & \frac{1}{4}\|\nabla^{4}
u\|_{L^{2}}^{2}+\frac{1}{4}\|\Lambda^{5} d\|_{L^{2}}^{2} +C_{0}(
\|\Delta d\|_{L^{2}}^{14}+\|\Delta d\|_{L^{2}}^{\frac{38}{7}}).
\end{align}
%--------------------(eq2.22)------------------------------------
Here we have used the following Gagliardo-Nirenberg inequalities:
\begin{align*}
&\|\Lambda^{4} d\|_{L^{3}}\leq C\|\Delta
d\|_{L^{2}}^{\frac{1}{6}}\|\Lambda^{5} d\|_{L^{2}}^{\frac{5}{6}};\\
& \|\Lambda^{2} d\|_{L^{4}}\leq C\|\Delta
d\|_{L^{2}}^{\frac{3}{4}}\|\Lambda^{5}d\|_{L^{2}}^{\frac{1}{4}};\\
&\|\Lambda^{3}d\|_{L^{4}}\leq C\|\Delta
d\|_{L^{2}}^{\frac{5}{6}}\|\Lambda^{5}d\|_{L^{2}}^{\frac{1}{6}}.
\end{align*}
Inserting \eqref{eq2.21} and \eqref{eq2.22} into \eqref{eq2.20}, one
gets
\begin{align}\label{eq2.23}
\frac{d}{dt}\|\Lambda^{3} u\|_{L^{2}}^{2}+&\|\Lambda^{4}
u\|_{L^{2}}^{2}\leq
\frac{1}{4}\|\Lambda^{5}d\|_{L^{2}}^{2}+C(\|\omega\|_{L^{2}}^{\frac{13}{2}}\|\Lambda^{3}
u\|_{L^{2}}^{\frac{3}{2}}+\|\Delta d\|_{L^{2}}^{14}+\|\Delta
d\|_{L^{2}}^{\frac{38}{7}})\nonumber\\
\leq& \frac{1}{4}\|\Lambda^{5}d\|_{L^{2}}^{2}+C_{0}C(\varepsilon)
(1+y(t))^{\frac{39C\varepsilon}{8}}\|\Lambda^{3}u\|_{L^{2}}^{\frac{3}{2}}+C_{0}C(\varepsilon)(1+y(t))^{\frac{21C\varepsilon}{2}}.
\end{align}
%--------------------(eq2.23)------------------------------------

Taking $\Lambda^{4}$ on \eqref{eq1.2}, multiplying $\Lambda^{4} d$
and integrating over $\mathbb{R}^{3}$, one obtains
\begin{align}\label{eq2.24}
\frac{1}{2}\frac{d}{dt}\|\Lambda^{4} d\|_{L^{2}}^{2}+\|\Lambda^{5}
d\|_{L^{2}}^{2}=-\int_{\mathbb{R}^{3}}\Lambda^{4}(u\cdot\nabla
d)\cdot\Lambda^{4}d\text{d}x+\int_{\mathbb{R}^{3}}\Lambda^{4}(|\nabla
d|^{2}d)\cdot\Lambda^{4}d\text{d}x:=I_{7}+I_{8}.
\end{align}
%--------------------(eq2.24)------------------------------------
Similar as the estimate of $I_{5}$, we have
\begin{align}\label{eq2.25}
I_{7}=&-\int_{\mathbb{R}^{3}}[\Lambda^{4} (u\cdot\nabla
d)-u\cdot\nabla\Lambda^{4} d]\cdot \Lambda^{4}d\text{d}x\nonumber\\
\leq&  C\|\Lambda^{4} (u\cdot\nabla d)-u\cdot\nabla\Lambda^{4}
d\|_{L^{\frac{3}{2}}}\|\Lambda^{4} d\|_{L^{3}}\nonumber\\
\leq& C\|\nabla d\|_{L^{6}}\|\Lambda^{4} u\|_{L^{2}}\|\Lambda^{4}
d\|_{L^{3}}+C\|\nabla u\|_{L^{6}}\|\Lambda^{4} d\|_{L^{2}}\|\Lambda^{4} d\|_{L^{3}}\nonumber\\
\leq& \frac{1}{4}\|\Lambda^{4} u\|_{L^{2}}^{2} +C\|\Delta
d\|_{L^{2}}^{2}\|\Lambda^{4} d\|_{L^{3}}^{2}+\|\Delta u\|_{L^{2}}\|\Lambda^{4} d\|_{L^{2}}\|\Lambda^{4} d\|_{L^{3}} \nonumber\\
\leq& \frac{1}{4}\|\Lambda^{4} u\|_{L^{2}}^{2} +C\|\Delta
d\|_{L^{2}}^{2}\|\Delta d\|_{L^{2}}^{\frac{1}{3}}\|\Lambda^{5}
d\|_{L^{2}}^{\frac{5}{3}}
+\|\nabla \omega \|_{L^{2}}\|\Delta d\|_{L^{2}}^{\frac{1}{2}}\|\Lambda^{5} d\|_{L^{2}}^{\frac{3}{2}}\nonumber\\
\leq & \frac{1}{4}\|\Lambda^{4} u\|_{L^{2}}^{2}+\frac{1}{4}
\|\Lambda^{5} d\|_{L^{2}}^{2} +C_{0}(\|\Delta
d\|_{L^{2}}^{14}+\|\nabla \omega\|_{L^{2}}^{14}+\|\Delta
d\|_{L^{2}}^{\frac{14}{5}}),
\end{align}
%--------------------(eq2.25)------------------------------------
where we have used the  Gagliardo-Nirenberg inequality:
\begin{align*}
\|\Lambda^{4} d\|_{L^{2}} \leq C\|\Delta
d\|_{L^{2}}^{\frac{1}{3}}\|\Lambda^{5}
d\|_{L^{2}}^{\frac{2}{3}}\text{ and }
 \|\Lambda^{4} d\|_{L^{3}} \leq
C\|\Delta d\|_{L^{2}}^{\frac{1}{6}}\|\Lambda^{5}
d\|_{L^{2}}^{\frac{5}{6}}.
\end{align*}
To estimate $I_{8}$, by using the Leibniz's rule, the fact $|d|=1$,
the H\"{o}lder's inequality               %, the interpolation inequalities
 and the
Young inequality, one obtains
\begin{align}\label{eq2.26}
I_{8}=&\int_{\mathbb{R}^{3}}\Lambda^{4}(|\nabla d|^{2}d)\cdot
\Lambda^{4}d\text{d}x=-\int_{\mathbb{R}^{3}}\Lambda^{3}(|\nabla
d|^{2}d)\cdot \Lambda^{5}d\text{d}x\nonumber\\
=&-\!\int_{\mathbb{R}^{3}}\!\big[\Lambda^{3}(|\nabla d|^{2})
d\cdot\Lambda^{5}d+3\Lambda^{2}(|\nabla d|^{2})\Lambda d
\cdot\Lambda^{5} d +3\Lambda (|\nabla d|^{2})\Lambda^{2} d \cdot
\Lambda^{5} d +|\nabla d|^{2}\Lambda^{3}d\cdot\Lambda^{5}
d\big]\text{d}x\nonumber\\
\leq & C\|\Lambda^{5} d\|_{L^{2}}(\|\nabla d\|_{L^{6}}\|\Lambda^{4}
d\|_{L^{3}}+\|\Lambda^{2}
d\|_{L^{4}}\|\Lambda^{3}d\|_{L^{4}}+\|\nabla
d\|_{L^{6}}^{2}\|\Lambda^{3} d\|_{L^{6}}+\|\nabla
d\|_{L^{6}}\|\Lambda^{2} d\|_{L^{6}}^{2})\nonumber\\
\leq & C\|\Lambda^{5} d\|_{L^{2}}(\|\Delta
d\|_{L^{2}}^{\frac{7}{6}}\|\Lambda^{5}d\|_{L^{2}}^{\frac{5}{6}} +\|
\Delta d\|_{L^{2}}^{\frac{7}{3}}\|\Lambda^{5}
d\|_{L^{2}}^{\frac{2}{3}})\nonumber\\
\leq &  \frac{1}{4}\|\Lambda^{5} d\|_{L^{2}}^{2}+C_{0}\|\Delta
d\|_{L^{2}}^{14}.
\end{align}
%--------------------(eq2.26)------------------------------------
Here we have used the following Gagliardo-Nirenberg inequalities:
\begin{align*}
&\|\Lambda^{4} d\|_{L^{3}}\leq C\|\Delta
d\|_{L^{2}}^{\frac{1}{6}}\|\Lambda^{5} d\|_{L^{2}}^{\frac{5}{6}};\\
& \|\Lambda^{2} d\|_{L^{4}}\leq C\|\Delta
d\|_{L^{2}}^{\frac{3}{4}}\|\Lambda^{5}d\|_{L^{2}}^{\frac{1}{4}};\\
&\|\Lambda^{3}d\|_{L^{6}}\leq C\|\Delta
d\|_{L^{2}}^{\frac{1}{3}}\|\Lambda^{5}d\|_{L^{2}}^{\frac{2}{3}};\\
&\|\Lambda^{2} d\|_{L^{6}}\leq C \|\Delta
d\|_{L^{2}}^{\frac{2}{3}}\|\Lambda^{5} d\|_{L^{2}}^{\frac{1}{3}}.
\end{align*}
Inserting \eqref{eq2.25} and \eqref{eq2.26} into \eqref{eq2.24}, one
gets
\begin{align}\label{eq2.27}
\frac{d}{dt}\|\Lambda^{4} d\|_{L^{2}}^{2}+\|\Lambda^{5}
d\|_{L^{2}}^{2}\leq& \frac{1}{2} \|\Lambda^{4}
d\|_{L^{2}}^{2}+C_{0}(\|\Delta d\|_{L^{2}}^{14}+\|\nabla \omega\|_{L^{2}}^{14}+\|\Delta d\|_{L^{2}}^{\frac{14}{5}})\nonumber\\
\leq&  \frac{1}{2} \|\Lambda^{4}
d\|_{L^{2}}^{2}+C_{0}C(\varepsilon)(1+y(t))^{\frac{21C\varepsilon}{2}}.
\end{align}
%--------------------(eq2.27)------------------------------------
Combining \eqref{eq2.23} and \eqref{eq2.27} together, one obtains
\begin{align*}
&\frac{d}{dt}(\|\Lambda^{3} u\|_{L^{2}}^{2}+\|\Lambda^{4}
d\|_{L^{2}}^{2})+\frac{1}{2}(\|\Lambda^{4}
u\|_{L^{2}}^{2}+\|\Lambda^{5} d\|_{L^{2}}^{2}) \nonumber\\
\leq& C_{0}C(\varepsilon)
(e+y(t))^{\frac{39C\varepsilon}{8}}\|\Lambda^{3}u\|_{L^{2}}^{\frac{3}{2}}+C_{0}C(\varepsilon)(e+y(t))^{\frac{21C\varepsilon}{2}}\nonumber\\
\leq& C_{0}C(\varepsilon)
(e+y(t))^{\frac{21C\varepsilon}{2}}(1+\|\Lambda^{3}u\|_{L^{2}}^{2}+\|\Lambda^{4}
d\|_{L^{2}}^{2})^{\frac{3}{4}}\nonumber\\
\leq&  C_{0}C(\varepsilon)
(e+y(t))^{\frac{3}{4}+\frac{21C\varepsilon}{2}}.
\end{align*}
Hence
\begin{align}\label{eq2.28}
\frac{d}{dt}y(t)\leq  C_{0}C(\varepsilon)
(e+y(t))^{\frac{3}{4}+\frac{21C\varepsilon}{2}}.
\end{align}
%--------------------(eq2.28)------------------------------------
By selecting $\varepsilon $ sufficiently small such that
$\frac{3}{4}+\frac{21C\varepsilon}{2}< 1$, then applying  Gronwall's
inequality to the above inequality \eqref{eq2.28}, we get the
boundness of $y(t)$ on $[T,T_{*}]$, i.e., the estimate \eqref{eq2.4}
is proved under the condition \eqref{eq2.3}. This completes the
proof of \eqref{eq1.6} under the assumptions of Theorem \ref{thm1.1}
.

\textit{Case II:} Next, we prove  \eqref{eq1.7} under the
assumptions of Theorem \ref{thm1.1}. Multiplying \eqref{eq1.1} by
$-\Delta u$, integrating over $\mathbb{R}^{3}$,  we get
\begin{align}\label{eq2.29}
\frac{1}{2}\frac{d}{dt}\|\nabla u (\cdot,t)\|_{L^{2}}^{2}
 +\|\Delta u (\cdot,t)\|_{L^{2}}^{2}=&-\int_{\mathbb{R}^{3}}(\nabla
 u\cdot\nabla)u\cdot\nabla u\text{d}x+\int_{\mathbb{R}^{3}}(\nabla
 d\cdot\Delta d)\cdot\Delta u\text{d}x\nonumber\\
=&-\int_{\mathbb{R}^{3}}(\nabla
 u\cdot\nabla)u\cdot\nabla u\text{d}x+\sum_{i,k=1}^{3}\int_{\mathbb{R}^{3}}\partial_{i}
 d_{k}\Delta d_{k}\Delta u_{i}\text{d}x.
\end{align}
%--------------------(eq2.29)------------------------------------
%Notice that the $i$th (i=1,2,3) component of $u$ satisfies
%\begin{align}\label{eq2.28}
%\partial_{t}u_{i}+(u\cdot\nabla)u_{i}-\Delta
%u_{i}+\partial_{i}P=-\sum_{j=1}^{3}\partial_{j}\left(\sum_{k=1}^{3}\partial_{i}d_{k}\partial_{j}d_{k}\right).
%\end{align}
%%--------------------(eq2.28)------------------------------------
%Multiplying (2.7) by $-\Delta u_{i}$, summing over $i$ and
%integrating by parts,  we get
%\begin{align}
%&\frac{1}{2}\frac{d}{dt}\|\nabla u (\cdot,t)\|_{L^{2}}^{2}
% +\|\Delta u (\cdot,t)\|_{L^{2}}^{2} \nonumber\\
%=&\sum_{i,j=1}^{3}\int_{\mathbb{R}^{3}}u_{j}\partial_{j} u_{i}
%\Delta
%u_{i}\text{d}x+\sum_{i,k=1}^{3}\int_{\mathbb{R}^{3}}\partial_{i}
%d_{k}\Delta d_{k}\Delta u_{i}\text{d}x\nonumber\\
%=&-\sum_{i,j=1}^{3}\int_{\mathbb{R}^{3}}\nabla (u_{j}\partial_{j}
%u_{i}) \nabla
%u_{i}\text{d}x+\sum_{i,k=1}^{3}\int_{\mathbb{R}^{3}}\partial_{i}
%d_{k}\Delta d_{k}\Delta u_{i}\text{d}x\nonumber\\
% =&-\sum_{i,j=1}^{3}\int_{\mathbb{R}^{3}}\nabla u_{j}\partial_{j}
%u_{i} \nabla
%u_{i}\text{d}x+\sum_{i,k=1}^{3}\int_{\mathbb{R}^{3}}\partial_{i}
%d_{k}\Delta d_{k}\Delta u_{i}\text{d}x,
%\end{align}
%%----------(2.8)------------
%where in the last equality, we used that fact that $\text{div } u=0$
%implies that
%$\sum_{i,j=1}^{3}\int_{\mathbb{R}^{3}}u_{j}\partial_{j}\nabla
%u_{i}\nabla
%u_{i}\text{d}x=\sum_{i,j=1}^{3}\int_{\mathbb{R}^{3}}u_{j}\partial_{j}(\frac{|\nabla
%u_{i}|^{2}}{2})\text{d}x=0$.
Applying $\Delta$ on \eqref{eq1.2}, multiplying it with $\Delta d$,
and using \eqref{eq1.3}, we obtain
\begin{align}\label{eq2.30}
\frac{1}{2}\frac{d}{dt}\|\Delta d (\cdot,t) \|_{L^{2}}^{2}&
+\|\nabla \Delta d\|_{L^{2}}^{2}
=-\!\int_{\mathbb{R}^{3}}\Delta(u\cdot\nabla d)\cdot\Delta
d\text{d}x+ \!\int_{\mathbb{R}^{3}}\Delta(|\nabla
d|^{2}d)\cdot\Delta d\text{d}x\nonumber\\
=& -\int_{\mathbb{R}^{3}}(\Delta u\cdot\nabla d)\cdot\Delta
d\text{d}x-2\int_{\mathbb{R}^{3}}(\nabla u\cdot \nabla)\nabla
d\cdot\Delta d\text{d}x+\int_{\mathbb{R}^{3}}\Delta (|\nabla
d|^{2}d)\cdot\Delta d\text{d}x\nonumber\\
=&\!-\!\!\sum_{i,k=1}^{3}\!\!\int_{\mathbb{R}^{3}}\!\!\partial_{i}
 d_{k}\Delta d_{k}\Delta u_{i}\text{d}x\!-\!2\!\!\int_{\mathbb{R}^{3}}\!\!(\nabla u\cdot \nabla)\nabla
d\cdot\Delta d\text{d}x+\!\int_{\mathbb{R}^{3}}\!\Delta (|\nabla
d|^{2}d)\cdot\Delta d\text{d}x.
\end{align}
%--------------------(eq2.30)------------------------------------
Combining \eqref{eq2.29} and \eqref{eq2.30} together, we have
\begin{align}\label{eq2.31}
&\frac{1}{2}\frac{d}{dt}\left(\|\nabla
u(\cdot,t)\|_{L^{2}}^{2}+\|\Delta d
(\cdot,t)\|_{L^{2}}^{2}\right)+(\|\Delta
u\|_{L^{2}}^{2}+\|\nabla \Delta d\|_{L^{2}}^{2})\nonumber\\
=&-\int_{\mathbb{R}^{3}}(\nabla
 u\cdot\nabla)u\cdot\nabla u\text{d}x-2\!\!\int_{\mathbb{R}^{3}}\!\!(\nabla u\cdot \nabla)\nabla
d\cdot\Delta d\text{d}x+\int_{\mathbb{R}^{3}}\Delta (|\nabla d|^{2}d)\Delta d\text{d}x\nonumber\\
=:&J_{1}+J_{2}+J_{3}.
\end{align}
%--------------------(eq2.31)------------------------------------
For $J_{1}$ and $J_{2}$, by the H\"{o}lder's inequality and
\eqref{eq2.1}, we have
\begin{align}
    \label{eq2.32}
&J_{1}\leq\int_{\mathbb{R}^{3}} |\nabla u||\nabla u| |\nabla u
|\text{d}x
\leq C\|\nabla u\|_{L^{4}}^{2}\|\nabla u\|_{L^{2}}\nonumber\\
\leq& C\|\Delta u\|_{L^{2}}\|\nabla
u\|_{\dot{B}^{-1}_{\infty,\infty}}\|\nabla u\|_{L^{2}}\nonumber\\
\leq& \frac{1}{2} \|\Delta u\|_{L^{2}}^{2}+C \|\nabla
u\|_{\dot{B}^{-1}_{\infty,\infty}}^{2}\|\nabla u\|_{L^{2}}^{2};\\
%--------------------(eq2.32)------------------------------------
    \label{eq2.33}
&J_{2}\leq C\|\nabla u\|_{L^{2}}\|\nabla^{2}d\|_{L^{4}}\|\Delta
d\|_{L^{4}}\leq C\|\nabla u\|_{L^{2}}\|\nabla
d\|_{L^{\infty}}\|\nabla\Delta d\|_{L^{2}}\nonumber\\
\leq & \frac{1}{4}\|\nabla \Delta d\|_{L^{2}}^{2}+C\|\nabla
d\|_{L^{\infty}}^{2}\|\nabla u\|_{L^{2}}^{2},
\end{align}
%--------------------(eq2.33)------------------------------------
where we have used  the Gagliardo-Nirenberg inequalities
\begin{align*}
\|\nabla^{2} d\|_{L^{4}}\leq C\|\nabla
d\|_{L^{\infty}}^{\frac{1}{2}}\|\nabla \Delta
d\|_{L^{2}}^{\frac{1}{2}}\text{ and } \|\Delta d\|_{L^{4}}\leq
C\|\nabla d\|_{L^{\infty}}^{\frac{1}{2}}\|\Delta
d\|_{L^{2}}^{\frac{1}{2}}.
\end{align*}
For $J_{3}$, similar as the estimate of $I_{4}$, we have
\begin{align}\label{eq2.34}
J_{3}%=&\int_{\mathbb{R}^{3}}(\Delta (|\nabla d|^{2}) d\cdot\Delta
%d+2\nabla (|\nabla d|^{2})\nabla d\cdot\Delta d+|\nabla
%d|^{2}|\Delta d|^{2})\text{d}x\nonumber\\
%\leq& \int_{\mathbb{R}^{3}}(2|\nabla d||\nabla\Delta d||\Delta
%d|+2|\nabla^{2} d|^{2}|\Delta d|+4|\nabla d|^{2}|\nabla^{2}
%d||\Delta d|+|\nabla d|^{2}|\Delta d|^{2})\text{d}x\nonumber\\
%\leq& C (\|\nabla \Delta d\|_{L^{2}}\|\nabla
%d\|_{L^{\infty}}\|\Delta d\|_{L^{2}} +\|\nabla^{2}
%d\|_{L^{4}}^{2}\|\Delta d\|_{L^{2}}+\|\nabla
%d\|_{L^{\infty}}\|\nabla^{2}d\|_{L^{2}}\|\Delta d\|_{L^{2}}+\|\nabla
%d\|_{L^{\infty}}^{2}\|\Delta d\|_{L^{2}}^{2})\nonumber\\
%\leq &C (\|\nabla \Delta d\|_{L^{2}}\|\nabla
%d\|_{L^{\infty}}\|\Delta d\|_{L^{2}} +\|\nabla
%d\|_{L^{\infty}}^{2}\|\Delta d\|_{L^{2}}^{2})\nonumber\\
\leq& \frac{1}{4}\|\nabla\Delta d\|_{L^{2}}^{2}+C\|\nabla
d\|_{L^{\infty}}^{2}\|\Delta d\|_{L^{2}}^{2}.
\end{align}
%--------------------(eq2.34)------------------------------------
Substituting the above estimates \eqref{eq2.32}--\eqref{eq2.34} into
\eqref{eq2.31}, we obtain
\begin{align}\label{eq2.35}
&\frac{d}{dt}\left(\|\nabla u(\cdot,t)\|_{L^{2}}^{2}+\|\Delta
d(\cdot,t)\|_{L^{2}}^{2}\right)+(\|\Delta
u\|_{L^{2}}^{2}+\|\nabla \Delta d\|_{L^{2}}^{2})\nonumber\\
\leq& C\|\nabla u\|_{\dot{B}^{-1}_{\infty,\infty}}^{2}\|\nabla
u\|_{L^2}^{2}+C\|\nabla d\|_{L^{\infty}}^{2}(\|\nabla
u\|_{L^{2}}^{2}+\|\Delta
d\|_{L^2}^{2})\nonumber\\
\leq&C(\|\nabla u\|_{\dot{B}^{-1}_{\infty,\infty}}^{2}+\|\nabla
d\|_{L^{\infty}}^{2})(\|\nabla u\|_{L^2}^{2}+\|\Delta
d\|_{L^2}^{2})\nonumber\\
\leq& C(1+\|\nabla u\|_{\dot{B}^{-1}_{\infty,\infty}}^{2}+\|\nabla
d\|_{\dot{B}^{0}_{\infty,\infty}}^{2}\text{ln}(e+\|\Lambda^{3}\nabla
d\|_{L^{2}}))(\|\nabla u\|_{L^2}^{2}+\|\Delta
d\|_{L^2}^{2})\nonumber\\
\leq & C\left(1+\frac{\|\nabla
u\|_{\dot{B}^{-1}_{\infty,\infty}}^{2}+\|\nabla
d\|_{\dot{B}^{0}_{\infty,\infty}}^{2}}{\sqrt{1+\text{ln}(e+\|\nabla
u\|_{\dot{B}^{-1}_{\infty,\infty}} +\|\nabla
d\|_{\dot{B}^{0}_{\infty,\infty}})}}\right)
\sqrt{1+\text{ln}(e+\|\nabla u\|_{\dot{B}^{-1}_{\infty,\infty}}
+\|\nabla
d\|_{\dot{B}^{0}_{\infty,\infty}})}\nonumber\\
&\times\text{ln}(e+\|\Lambda^{3} u\|_{L^{2}}+\|\Lambda^{4}
d\|_{L^{2}}))(\|\nabla u\|_{L^2}^{2}+\|\Delta
d\|_{L^2}^{2})\nonumber\\
\leq & C\left(1+\frac{\|\nabla
u\|_{\dot{B}^{-1}_{\infty,\infty}}^{2}+\|\nabla
d\|_{\dot{B}^{0}_{\infty,\infty}}^{2}}{\sqrt{1+\text{ln}(e+\|\nabla
u\|_{\dot{B}^{-1}_{\infty,\infty}} +\|\nabla
d\|_{\dot{B}^{0}_{\infty,\infty}})}}\right) \sqrt{1+\text{ln}(e+\|
u\|_{\dot{B}^{0}_{\infty,\infty}} +\|\nabla
d\|_{\dot{B}^{0}_{\infty,\infty}})}\nonumber\\
&\times\text{ln}(e+\|\Lambda^{3} u\|_{L^{2}}+\|\Lambda^{4}
d\|_{L^{2}}))(\|\nabla u\|_{L^2}^{2}+\|\Delta
d\|_{L^2}^{2})\nonumber\\
\leq& C \!\left(\!1+\!\frac{\|\nabla
u\|_{\dot{B}^{-1}_{\infty,\infty}}^{2}+\|\nabla
d\|_{\dot{B}^{0}_{\infty,\infty}}^{2}}{\sqrt{1+\text{ln}(e+\|\nabla
u\|_{\dot{B}^{-1}_{\infty,\infty}} +\|\nabla
d\|_{\dot{B}^{0}_{\infty,\infty}})}}\right)
(1+\text{ln}(e+y(t)))^{\frac{3}{2}}(\|\nabla u\|_{L^2}^{2}+\|\Delta
d\|_{L^2}^{2}),
\end{align}
%--------------------(eq2.35)------------------------------------
where $y(t)$ is defined as \eqref{eq2.15}. Here we have used the
Sobolev imbedding \eqref{eq2.8} again. Similarly, assume that the
condition \eqref{eq1.7} is not ture, one concludes that for any
small constant $\varepsilon>0$, there exists
$T=T(\varepsilon)\in(0,T_{*})$ such that
\begin{align}\label{eq2.36}
\int_{T}^{T_{*}}\frac{\|\nabla
u\|_{\dot{B}^{-1}_{\infty,\infty}}^{2}+\|\nabla
d\|_{\dot{B}^{0}_{\infty,\infty}}^{2}}{\sqrt{1+\text{ln}(e+\|\nabla
u\|_{\dot{B}^{-1}_{\infty,\infty}} +\|\nabla
d\|_{\dot{B}^{0}_{\infty,\infty}})}}\text{d}t\leq \varepsilon.
\end{align}
%--------------------(eq2.36)------------------------------------
Applying Gronwall's inequality on \eqref{eq2.36} for the interval
$[T,t)$, one obtains
\begin{align}\label{eq2.37}
&\|\nabla u (\cdot,t)\|_{L^{2}}^{2}+\|\Delta d
(\cdot,t)\|_{L^{2}}^{2}\nonumber\\
  \leq& (\|\nabla u
(\!\cdot,\!T)\|_{\!L^{\!2}}^{\!2}\!\!+\!\|\Delta d
(\!\cdot,\!T)\|_{\!L^{\!2}}^{\!2})\exp\!\!\left\{\!\!\int_{T}^{t}\!\!C\text{d}\tau\!+C\!\!\int_{T}^{t}\!\!
\frac{\|\nabla u\|_{\dot{B}^{-1}_{\infty,\infty}}^{2}+\|\nabla
d\|_{\dot{B}^{0}_{\infty,\infty}}^{2}}{\!\sqrt{1\!+\!\text{ln}(e\!+\!\|\nabla
u\|_{\!\dot{B}^{-1}_{\!\infty,\!\infty}}
\!\!+\!\|\nabla d\|_{\!\dot{B}^{0}_{\!\infty,\!\infty}}\!)}}\text{d}\tau\text{ln}^{\frac{3}{2}}(e\!+y(t))\!\!\right\}\nonumber\\
 \leq& (\|\nabla u
(\cdot,T)\|_{L^{2}}^{2}+\|\Delta d
(\cdot,T)\|_{L^{2}}^{2})\exp\left\{C[(T_{*}-T)+\varepsilon\text{ln}^{\frac{3}{2}}(e+y(t))]\right\}\nonumber\\
\leq &
C(\varepsilon)\exp(C\varepsilon(1+\text{ln}^{\frac{3}{2}}(e+y(t))))
\leq C(\varepsilon) (e+y(t))^{\frac{3C\varepsilon}{2}}.
\end{align}
%--------------------(eq2.37)------------------------------------
Similar as the proofs in \textit{Case I}, by using the above
estimate \eqref{eq2.37}, we find that \eqref{eq2.28} still remains
valid. Thus, we deduce that the estimate \eqref{eq2.4} still holds
provided that the $\varepsilon$ in \eqref{eq2.36} sufficiently
small, i.e., \eqref{eq2.4} holds under the condition \eqref{eq1.7}
is not true, which implies that $[0,T_{*})$ is not a maximal
interval of existence of smooth solutions and leads to the
contradiction. This completes the proof of
\eqref{eq1.7}.\hfill$\Box$

\section{Proof of Theorem \ref{cor1.3}}

Similar as the proof of Theorem \ref{thm1.1}, we prove Theorem
\ref{cor1.3} by contradiction. Assume that \eqref{eq1.8} is not
true, then there exists $0<M<\infty$ such that
\begin{align}\label{eq3.1}
\int_{0}^{T_{*}}\frac{\|\nabla
d\|_{\dot{B}^{0}_{\infty,\infty}}^{2}}{\sqrt{1+\text{ln}(e +\|\nabla
d\|_{\dot{B}^{0}_{\infty,\infty}})}}\text{d}t\leq M.
\end{align}
%--------------------(eq3.1)------------------------------------
In what follows, we will give some a priori estimates to show that
\begin{align}\label{eq3.2}
\lim_{t\rightarrow T_{*}^{-}}(\|\nabla^{2}
u(t,\cdot)\|_{L^2}^{2}+\|\nabla^{3} d(t,\cdot)\|_{L^2}^{2})\leq C
\end{align}
%--------------------(eq3.2)------------------------------------
for some positive constant $C$ depends only on $u_{0}, d_{0}, T_{*}$
and $M$.

We first taking $\nabla\times$ on \eqref{eq1.1}, it follows from
dimension $n=2$ that
\begin{align*}
\omega_{t}-\Delta \omega+ u\cdot\nabla \omega=-\nabla\times(\Delta
d\cdot\nabla d).
\end{align*}
Multiplying above equality with $\omega$ and integrating over
$\mathbb{R}^{2}$, we obtain
\begin{align}\label{eq3.3}
\frac{1}{2}\frac{d}{dt}\|\omega(\cdot,t)\|_{L^{2}}^{2}+&\|\nabla
\omega\|_{L^{2}}^{2}=-\int_{\mathbb{R}^{2}}\nabla\times(\Delta
d\cdot\nabla d)\cdot\omega\text{d}x\nonumber\\
=& \int_{\mathbb{R}^{2}}(\Delta d\cdot\nabla
d)\cdot\nabla\times\omega\text{d}x\leq C\|\nabla
d\|_{L^{\infty}}\|\Delta d\|_{L^{2}}\|\nabla \omega\|_{L^{2}}
\nonumber\\
\leq& \frac{1}{4}\|\nabla\omega\|_{L^{2}}^{2}+C\|\nabla
d\|_{L^{\infty}}^{2}\|\Delta d\|_{L^{2}}^{2}.
\end{align}
%--------------------(eq3.3)------------------------------------
Similar as the estimate of $\Delta d$ on the \textit{Case I} in
section 2, we have
\begin{align}\label{eq3.4}
\frac{1}{2}\frac{d}{dt}&\|\Delta d (\cdot,t)\|_{L^{2}}^{2}+\|\nabla
\Delta d\|_{L^{2}}^{2}=-\!\int_{\mathbb{R}^{2}}\Delta(u\cdot\nabla
d)\cdot\Delta d\text{d}x+ \!\int_{\mathbb{R}^{2}}\Delta(|\nabla
d|^{2}d)\cdot\Delta d\text{d}x\nonumber\\
\leq& \int_{\mathbb{R}^{2}}(|\Delta u||\nabla d||\Delta d|+2|\nabla
u||\nabla^{2}d||\Delta d|+\Delta(|\nabla d|^{2})d\cdot\Delta
d+2\nabla(|\nabla d|^{2})\nabla d\cdot\Delta d+|\nabla d|^{2}|\Delta
d|^{2})\text{d}x\nonumber\\
\leq & C(\|\Delta u\|_{L^{2}}\|\nabla d\|_{L^{\infty}}\|\Delta
d\|_{L^{2}}+\|\nabla u\|_{L^{2}}\|\nabla^{2} d\|_{L^{4}}\|\Delta
d\|_{L^{4}}\nonumber\\
&+\|\nabla d\|_{L^{\infty}}^{2}\|\Delta
d\|_{L^{2}}\|\nabla^{2}d\|_{L^{2}}+\|\nabla
d\|_{L^{\infty}}\|\nabla^{2}d\|_{L^{2}}\|\nabla \Delta
d\|_{L^{2}}+\|\nabla d\|_{L^{\infty}}^{2}\|\Delta
d\|_{L^{2}}^{2})\nonumber\\
\leq& C (\|\nabla \omega\|_{L^{2}}\|\nabla d\|_{L^{\infty}}\|\Delta
d\|_{L^{2}}+\|\omega\|_{L^{2}}\|\nabla d\|_{L^{\infty}}\|\nabla
\Delta d\|_{L^{2}}+\|\nabla d\|_{L^{\infty}}^{2}\|\Delta
d\|_{L^{2}}^{2})\nonumber\\
\leq &\frac{1}{4}\|\nabla\omega\|_{L^{2}}^{2}+\frac{1}{2}\|\nabla
\Delta d\|_{L^{2}}^{2}+C\|\nabla
d\|_{L^{\infty}}^{2}(\|\omega\|_{L^{2}}^{2}+\|\Delta
d\|_{L^{2}}^{2}),
\end{align}
%--------------------(eq3.4)------------------------------------
where we have used the facts $|d|=1$, $\|\Delta u\|_{L^{2}}\leq
C\|\nabla \omega\|_{L^{2}}$,  $\|\nabla u\|_{L^{2}}\leq
C\|\omega\|_{L^{2}}$ and $\|\nabla^{2}d\|_{L^{2}}=\|\Delta
d\|_{L^{2}}$, and the Gagliardo-Nirenberg inequalities in
$\mathbb{R}^{2}$:
\begin{align*}
\|\nabla^{2} d\|_{L^{4}}\leq C\|\nabla
d\|_{L^{\infty}}^{\frac{1}{2}}\|\nabla \Delta
d\|_{L^{2}}^{\frac{1}{2}}\text{ and } \|\Delta d\|_{L^{4}}\leq
C\|\nabla d\|_{L^{\infty}}^{\frac{1}{2}}\|\Delta
d\|_{L^{2}}^{\frac{1}{2}}.
\end{align*}
Notice that due to \eqref{eq3.1}, one concludes that for any small
constant $\varepsilon>0$, there exists $T=T(\varepsilon)\in
(0,T_{*})$ such that
\begin{align*}
\int_{T}^{T_{*}}\frac{\|\nabla
d\|_{\dot{B}^{0}_{\infty,\infty}}^{2}}{\sqrt{1+\text{ln}(e+\|\nabla
d\|_{\dot{B}^{0}_{\infty,\infty}})}}\text{d}t\leq \varepsilon.
\end{align*}
For any $T<t\leq T_{*}$, let
\begin{align}\label{eq3.5}
z(t):=\sup_{T<\tau\leq t}(\|\Lambda^{2}u\|_{L^{2}}^{2}+\|\Lambda^{3}
d\|_{L^{2}}^{2}).
\end{align}
%--------------------(eq3.5)------------------------------------
Combining \eqref{eq3.3} and \eqref{eq3.4} together, and using the
inequality \eqref{eq2.2} with $s=3$, it follows that
\begin{align}\label{eq3.6}
\frac{d}{dt}&(\|\omega (\cdot,t)\|_{L^{2}}^{2}+\|\Delta
d\|_{L^{2}}^{2})+(\|\nabla \omega \|_{L^{2}}^{2}+\|\nabla\Delta
d\|_{L^{2}}^{2})\leq C\|\nabla d\|_{L^{\infty}}^{2}(\|\omega
\|_{L^{2}}^{2}+\|\Delta d \|_{L^{2}}^{2})\nonumber\\
\leq &C(1+\|\nabla
d\|_{\dot{B}^{0}_{\infty,\infty}}^{2}\text{ln}(e+\|\Lambda^{2}\nabla
d\|_{L^{2}}))(\|\omega \|_{L^{2}}^{2}+\|\Delta d
\|_{L^{2}}^{2})\nonumber\\
\leq & C(1+\frac{\|\nabla
d\|_{\dot{B}^{0}_{\infty,\infty}}^{2}}{\sqrt{1+\text{ln}(e+\|\nabla
d\|_{\dot{B}^{0}_{\infty,\infty}})}})\sqrt{1+\text{ln}(e+\|\nabla
d\|_{\dot{B}^{0}_{\infty,\infty}})}\text{ln}(e+\|\Lambda^{3}
d\|_{L^{2}}))(\|\omega \|_{L^{2}}^{2}+\|\Delta d
\|_{L^{2}}^{2})\nonumber\\
\leq & C(1+\frac{\|\nabla
d\|_{\dot{B}^{0}_{\infty,\infty}}^{2}}{\sqrt{1+\text{ln}(e+\|\nabla
d\|_{\dot{B}^{0}_{\infty,\infty}})}})(\text{ln}(e+z(t)))^{\frac{3}{2}}(\|\omega
\|_{L^{2}}^{2}+\|\Delta d \|_{L^{2}}^{2}),
\end{align}
%--------------------(eq3.6)------------------------------------
where we have used the Sobolev imbedding
\begin{align*}
H^{2}(\mathbb{R}^{2})\subseteq L^{\infty}(\mathbb{R}^{2})\subseteq
\dot{B}^{0}_{\infty,\infty}(\mathbb{R}^{2}).
\end{align*}
By Gronwall's inequality to \eqref{eq3.6} in the interval $[T,t)$,
we have
\begin{align}\label{eq3.7}
&\|\omega(\cdot,t)\|_{L^{2}}^{2}+\|\Delta d (\cdot,
t)\|_{L^{2}}^{2}\nonumber\\
 \leq& (\|\omega
(\!\cdot,\!T)\|_{\!L^{\!2}}^{\!2}\!\!+\|\Delta d
(\!\cdot,\!T)\|_{\!L^{\!2}}^{\!2})\exp\!\!
\left\{\int_{T}^{t}C\text{d}\tau+C\int_{T}^{t}\frac{\|\nabla
d\|_{\dot{B}^{0}_{\infty,\infty}}^{2}}{\sqrt{1+\text{ln}(e+\|\nabla
d\|_{\dot{B}^{0}_{\infty,\infty}})}}\text{d}\tau(\text{ln}(e+y(t)))^{\frac{3}{2}}
\right\}\nonumber\\
 \leq& (\|\omega
(\cdot,T)\|_{L^{2}}^{2}+\|\Delta d
(\cdot,T)\|_{L^{2}}^{2})\exp\left\{C[(T_{*}-T)+\varepsilon\text{ln}^{\frac{3}{2}}(e+z(t))]\right\}\nonumber\\
\leq &
C(\varepsilon)\exp(C\varepsilon(1+\text{ln}^{\frac{3}{2}}(e+z(t))))
\leq C(\varepsilon) (e+z(t))^{\frac{3C\varepsilon}{2}}.
\end{align}
%--------------------(eq3.1)------------------------------------

On the other hand, by Lin, Lin and Wang \cite{LLW}, we known that
when the dimension $n=2$, there holds the energy equality
\begin{align}\label{eq3.8}
\|u(\cdot,t)\|_{L^{2}}^{2}+\|\nabla
d(\cdot,t)\|_{L^{2}}^{2}+\int_{0}^{T_{*}}\int_{\mathbb{R}^{2}}(|\nabla
u|^{2}+|\Delta d+|\nabla
d|^{2}d|^{2})\text{d}x\text{d}t=\|u_{0}\|_{L^{2}}^{2}+\|\nabla
d_{0}\|_{L^{2}}^{2}.
\end{align}
%--------------------(eq3.8)------------------------------------

Now, we will estimate $z(t)$ defined by \eqref{eq3.5}. Let us recall
the following useful Gagliardo-Nirenberg inequalities in
$\mathbb{R}^{2}$:
\begin{align}\label{eq3.9}
&\|\nabla u\|_{L^{3}}\leq C\|\nabla
u\|_{L^{2}}^{\frac{5}{6}}\|\Lambda^{3} u\|_{L^{2}}^{\frac{1}{6}};\quad\quad %\nonumber\\
\|\nabla u\|_{L^{6}}\leq C\|
u\|_{L^{2}}^{\frac{4}{9}}\|\Lambda^{3} u\|_{L^{2}}^{\frac{5}{9}};\nonumber\\
&\|\Lambda^{2} u\|_{L^{3}}\leq C\|\nabla
u\|_{L^{2}}^{\frac{1}{3}}\|\Lambda^{3} u\|_{L^{2}}^{\frac{2}{3}};\quad\quad %\nonumber\\
\|\nabla d\|_{L^{6}}\leq C\|\nabla
d\|_{L^{2}}^{\frac{1}{3}}\|\Delta d\|_{L^{2}}^{\frac{2}{3}};\\
&\|\Lambda^{2} d\|_{L^{4}}\leq C\|\Delta
d\|_{L^{2}}^{\frac{3}{4}}\|\Lambda^{4} d\|_{L^{2}}^{\frac{1}{4}};\quad\quad %\nonumber\\
\|\Lambda^{3} d\|_{L^{3}}\leq C\|\Delta
d\|_{L^{2}}^{\frac{1}{3}}\|\Lambda^{4}
d\|_{L^{2}}^{\frac{2}{3}}.\nonumber
\end{align}
%--------------------(eq3.9)------------------------------------
Applying $\Lambda^{2}$ on \eqref{eq1.1}, multiplying $\Lambda^{2} u$
and integrating over $\mathbb{R}^{3}$, and using \eqref{eq2.18}, the
H\"{o}lder's inequality, \eqref{eq3.9} and the Young inequalities,
one obtains
\begin{align}\label{eq3.10}
&\frac{1}{2}\frac{d}{dt}\|\Lambda^{2}u(\cdot,t)\|_{L^{2}}^{2}+\|\Lambda^{3}u\|_{L^{2}}^{2}
=-\int_{\mathbb{R}^{2}}\Lambda^{2}(u\cdot\nabla
u)\cdot\Lambda^{2}u\text{d}x-\int_{\mathbb{R}^{2}}\Lambda^{2}(\Delta
d\cdot \nabla d)\cdot\Lambda^{2}u\text{d}x\nonumber\\
=&-\int_{\mathbb{R}^{2}}[\Lambda^{2}(u\cdot\nabla u)-u\cdot\nabla
\Lambda^{2} u]\cdot\Lambda^{2}
u\text{d}x+\int_{\mathbb{R}^{2}}\Lambda(\Delta
d\cdot \nabla d)\cdot\Lambda^{3}u\text{d}x\nonumber\\
\leq & C\|[\Lambda^{2}(u\cdot\nabla u)-u\cdot\nabla \Lambda^{2}
u]\|_{L^{\frac{3}{2}}}\|\Lambda^{2} u\|_{L^{3}}+C(\|\nabla
d\|_{L^{6}}\|\Lambda^{3}d\|_{L^{3}}\|\Lambda^{3} u\|_{L^{2}}
+\|\Lambda^{3}u\|_{L^{2}}\|\Lambda^{2} d\|_{L^{4}}^{2})\nonumber\\
\leq& C (\|\Lambda^{2} u\|_{L^{3}}^{2}\|\nabla u\|_{L^{3}}+\|\nabla
d\|_{L^{6}}\|\Lambda^{3}d\|_{L^{3}}\|\Lambda^{3} u\|_{L^{2}}
+\|\Lambda^{3}u\|_{L^{2}}\|\Lambda^{2} d\|_{L^{4}}^{2})\nonumber\\
\leq &C (\|\nabla
u\|_{\!L^{\!2}}^{\frac{3}{2}}\|\Lambda^{3}u\|_{\!L^{\!2}}^{\frac{3}{2}}\!+\|\Lambda^{3}
u\|_{\!L^{\!2}}\|\nabla d\|_{\!L^{\!2}}^{\frac{1}{3}}\|\Delta
d\|_{\!L^{2}}\|\Lambda^{4}
d\|_{\!L^{\!2}}^{\frac{2}{3}}\!+\|\Lambda^{3}
u\|_{\!L^{\!2}}\|\Delta d\|_{\!L^{\!2}}^{\frac{3}{2}}\|\Lambda^{4}
d\|_{\!L^{\!2}}^{\frac{1}{2}})\nonumber\\
\leq & \frac{1}{4}\|\Lambda^{3}
u\|_{L^{2}}^{2}+\frac{1}{4}\|\Lambda^{4} d\|_{L^{2}}^{2}+C(\|\nabla
u\|_{L^{2}}^{6}+ \|\nabla d\|_{L^{2}}^{2}\|\Delta
d\|_{L^{2}}^{6}+\|\Delta d\|_{L^{2}}^{6})\nonumber\\
\leq & \frac{1}{4}\|\Lambda^{3}
u\|_{L^{2}}^{2}+\frac{1}{4}\|\Lambda^{4}
d\|_{L^{2}}^{2}+C(1+\|\nabla u\|_{L^{2}}^{6}+ \|\Delta
d\|_{L^{2}}^{6})\nonumber\\
\leq & \frac{1}{4}\|\Lambda^{3}
u\|_{L^{2}}^{2}+\frac{1}{4}\|\Lambda^{4}
d\|_{L^{2}}^{2}+C(1+\|\omega \|_{L^{2}}^{6}+ \|\Delta
d\|_{L^{2}}^{6}),
\end{align}
%--------------------(eq3.10)------------------------------------
where we have used the energy equality \eqref{eq3.8}.  Taking
$\Lambda^{3}$ on \eqref{eq1.2}, multiplying $\Lambda^{3} d$,
integrating over $\mathbb{R}^{2}$, and using \eqref{eq2.18}, the
H\"{o}lder's inequality, \eqref{eq3.9} and the Young inequalities,
one obtains
\begin{align}\label{eq3.11}
\frac{1}{2}&\frac{d}{dt}\|\Lambda^{3}d(\cdot,t)\|_{L^{2}}^{2}+\|\Lambda^{4}
d \|_{L^{2}}^{2}=-\int_{\mathbb{R}^{2}}\Lambda^{3}(u\cdot\nabla
d)\cdot
\Lambda^{3}d\text{d}x+\int_{\mathbb{R}^{2}}\Lambda^{3}(|\nabla
d|^{2}d)\cdot\Lambda^{3}d\text{d}x\nonumber\\
=& -\int_{\mathbb{R}^{2}}[\Lambda^{3}(u\cdot\nabla
d)-u\cdot\nabla\Lambda^{3}d]\cdot
\Lambda^{3}d\text{d}x-\int_{\mathbb{R}^{2}}\Lambda^{2}(|\nabla
d|^{2}d)\cdot\Lambda^{4}d\text{d}x\nonumber\\
=& -\!\!\int_{\mathbb{R}^{2}}\![\Lambda^{3}(u\cdot\nabla
d)-u\cdot\nabla\Lambda^{3}d]\cdot
\Lambda^{3}d\text{d}x\!-\!\!\int_{\mathbb{R}^{2}}\![\Lambda^{2}(|\nabla
d|^{2})d+\!2\Lambda(|\nabla d|^{2})\Lambda d
+\!|\nabla d|^{2}\Lambda^{2} d]\cdot\Lambda^{4}d\text{d}x\nonumber\\
\leq& C\|[\Lambda^{3}(u\cdot\nabla
d)-u\cdot\nabla\Lambda^{3}d]\|_{L^{\frac{3}{2}}}\|\Lambda^{3}d\|_{L^{3}}
+C(\|\Lambda^{3} d\|_{L^{3}}\|\nabla d\|_{L^{6}}+\|\Lambda^{2}
d\|_{L^{4}}^{2}\nonumber\\
&+\|\Lambda d\|_{L^{6}}^{2}\|\Lambda^{2} d\|_{L^{6}}+\|\nabla
d\|_{L^{6}}\|\Lambda^{3}d\|_{L^{3}})\|\Lambda^{4}
d\|_{L^{2}}\nonumber\\
\leq &C (\|\Lambda^{3} d\|_{L^{3}}\|\nabla d\|_{L^{6}}\|\Lambda^{3}
u\|_{L^{2}}+\|\Lambda^{3} d\|_{L^{3}}^{2}\|\nabla u\|_{L^{6}}+
\|\Lambda^{3} d\|_{L^{3}}\|\nabla d\|_{L^{6}}\|\Lambda^{4}
d\|_{L^{2}}\nonumber\\
&+\|\Lambda^{2} d\|_{L^{4}}^{2}\|\Lambda^{4} d\|_{L^{2}}+\|\Lambda
d\|_{L^{6}}^{2}\|\Lambda^{2} d\|_{L^{6}}\|\Lambda^{4}
d\|_{L^{2}})\nonumber\\
\leq & \frac{1}{8}\|\Lambda^{3}
u\|_{L^{2}}^{2}+\frac{1}{8}\|\Lambda^{4} d\|_{L^{2}}^{2} +C(\|\nabla
d\|_{L^{6}}^{2}\|\Lambda^{3} d\|_{L^{3}}^{2}+\|\nabla
u\|_{L^{6}}\|\Lambda^{3} d\|_{L^{3}}^{2}+\|\Lambda^{2}
d\|_{L^{4}}^{4}+\|\Lambda d\|_{L^{6}}^{4}\|\Lambda^{2}
d\|_{L^{6}}^{2})\nonumber\\
\leq & \frac{1}{8}\|\Lambda^{3}
u\|_{L^{2}}^{2}+\frac{1}{8}\|\Lambda^{4} d\|_{L^{2}}^{2} +C(\|\nabla
d\|_{L^{2}}^{\frac{2}{3}}\|\Delta d\|_{L^{2}}^{2}\|\Lambda^{4}
d\|_{L^{2}}^{\frac{4}{3}} +\|u\|_{L^{2}}^{\frac{4}{9}} \|\Lambda^{3}
u\|_{L^{2}}^{\frac{5}{9}} \|\Delta
d\|_{L^{2}}^{\frac{2}{3}}\|\Lambda^{4}
d\|_{L^{2}}^{\frac{4}{3}}\nonumber\\
& +\|\Delta d\|_{L^{2}}^{3}\|\Lambda^{4} d\|_{L^{2}}+\|\nabla
d\|_{L^{2}}^{\frac{4}{3}}\|\Delta d\|_{L^{2}}^{4}\|\Lambda^{4}
d\|_{L^{2}}^{\frac{2}{3}})\nonumber\\
\leq & \frac{1}{4}\|\Lambda^{3}
u\|_{L^{2}}^{2}+\frac{1}{4}\|\Lambda^{4} d\|_{L^{2}}^{2} +
C(\|\nabla d\|_{L^{2}}^{2}\|\Delta
d\|_{L^{2}}^{6}+\|u\|_{L^{2}}^{8}\|\Delta d\|_{L^{2}}^{12}+\|\Delta
d\|_{L^{2}}^{6}+\|\nabla d\|_{L^{2}}^{2}\|\Delta d\|_{L^{2}}^{6}
)\nonumber\\
\leq & \frac{1}{4}\|\Lambda^{3}
u\|_{L^{2}}^{2}+\frac{1}{4}\|\Lambda^{4} d\|_{L^{2}}^{2} +
C(\|\Delta d\|_{L^{2}}^{6}+\|\Delta
d\|_{L^{2}}^{12})\nonumber\\
\leq& \frac{1}{4}\|\Lambda^{3}
u\|_{L^{2}}^{2}+\frac{1}{4}\|\Lambda^{4} d\|_{L^{2}}^{2} +
C(1+\|\Delta d\|_{L^{2}}^{12}),
\end{align}
%--------------------(eq3.11)------------------------------------
where we have used the equality \eqref{eq3.8}. Combining
\eqref{eq3.10} and \eqref{eq3.11} together, and using estimate
\eqref{eq3.7}, we obtain
\begin{align}\label{eq3.12}
&\frac{d}{dt}(\|\Lambda^{2}u(\cdot,t)\|_{L^{2}}^{2}+|\Lambda^{3}d(\cdot,t)\|_{L^{2}}^{2})
+(\|\Lambda^{3}u\|_{L^{2}}^{2} +\|\Lambda^{4} d
\|_{L^{2}}^{2})\nonumber\\
\leq & C (1+\|\omega\|_{L^{2}}^{6}+\|\Delta d\|_{L^{2}}^{12}) \leq C
(1+\|\omega\|_{L^{2}}^{12}+\|\Delta d\|_{L^{2}}^{12})
\nonumber\\
\leq& C(\varepsilon) (1+z(t))^{9C\varepsilon}.
\end{align}
%--------------------(eq3.12)------------------------------------
By selecting $\varepsilon $ sufficiently small such that
${9C\varepsilon}< 1$, and applying  Gronwall's inequality to the
above inequality \eqref{eq3.12}, we get the boundness of $z(t)$ on
$[T,T_{*}]$, i.e., the estimate \eqref{eq3.2} is proved under the
assumption that \eqref{eq1.8} is not true, which implies that
$[0,T_{*})$ is not a maximal interval of existence of smooth
solutions and leads to the contradiction. This completes the proof
of Theorem \ref{cor1.3}. \hfill$\Box$

%--------------------(proof of thm1.4)------------------------------------
%\\
%\\
%\textbf{Acknowledgments}


\begin{thebibliography}{99}

\bibitem{BKM} J. Beale, T. Kato and A. Majda, Remarks on breakdown of smooth solutions for
the 3D Euler equations. Commun. Math. Phys., 94 (1984), 61--66.

\bibitem{CDY} K. Chang, W. Ding and R. Ye, Finite-time blow-up of
the heat flow of harmonic maps from surfaces, J. Differ. Geom.,
36(2) (1992), 507--515.

\bibitem{CS} Y. Chen and M. Struwe, Existence and partial regularity results for the heat
 flow of harmonic maps, Math. Z., 201 (1989), 83103.


\bibitem{ER} J. L. Ericksen, Hydrostatic theory of liquid crystal, Arch. Rational Mech. Anal.,
 9 (1962), 371--378.


%\bibitem{HKL} R. Hardt, D. Kinderleher and F. H. Lin, Existence and
%partial regularity of static liquid crystal configurations, Commun.
%Math. Phys., 105 (1986), 547--570.

\bibitem{JHW} J. Hineman and C. Wang, Well-posedness of nematic
liquid crystal flow in $L^{3}_{loc}(\mathbb{R}^{3})$,
arXiv:1208.5965v1 [math.AP] 29 Aug. 2012.


\bibitem{HW} X. Hu and D. Wang, Global Solution to the
Three-Dimensional Incompressible Flow of Liquid Crystals, Commun.
Math. Phys., 296 (2010), 861--880.

\bibitem{HW1} T. Huang and C. Wang, Blow up criterion for Nematic
Liquid Crystal Flows, Comm. Partial Differ. Equ., 37 (2012),
875--884.

\bibitem{YG} Y. Giga, Solutions for semilinear parabolic equations in $L^{p}$ and regularity of weak solutions of the
Navier--Stokes system, J. Diff. Equ., 61 (1986), 186--212.

\bibitem{GG11} Z. Guo and S. Gala, Remarks on logarithmical
regularity criteria for the Navier-Stokes equations, J, Math. Phys.,
52 (2011), 063503.

\bibitem{KP} T. Kato and G. Ponce, Commutator estimates and the
Euler and Navier--Stokes equations, Commu. Pure Appl. Math., 41
(1988), 891--907.

\bibitem{KOT} H. Kozono, T. Ogawa and Y. Taniuchi, The critical
Sobolev inequalities in Besov spaces and regularity criterion to
some semi-linear evolution equations, Math. Z., 242 (2002),
251--278.

\bibitem{KT} H. Kozono and Y. Taniuchi, Bilinear estimates in BMO
and the Navier--Stokes equations, Math. Z., 235 (2000), 173--194.

\bibitem{LE} F. Leslie, Theory of flow phenomenum in liquid crystals. In: The Theory of
Liquid Crystals, London-New York: Academic Press, 4 (1979), 1--81.

\bibitem{PG}  P.G. Lemari\'{e}-Rieusset, \textit{Recent Developments in the
 Navier-Stokes Problem}, Chapman and Hall/CRC, 2002.

\bibitem{XLW} X. Li and D. Wang, Global solution to the
incompressible flow of liquid crystal, J. Differ. Equ., 252 (2012),
745--767.

\bibitem{L} F. Lin, Nonlinear theory of defects in nematic liquid crystals; phase transition and flow phenomena,
 Comm. Pure. Appl. Math., 42 (1989), 789--814.

 \bibitem{LLW} F. Lin, J. Lin and C. Wang, Liquid Crystal flow in two dimensions, Arch.
 Rational Mech. Anal., 197 (2010), 297--336

\bibitem{LL1} F. Lin and C. Liu, Nonparabolic dissipative systems modeling the flow of liquid crystals,
Comm. Pure. Appl. Math.,  48  (1995), 501--537.

\bibitem{LL2} F. Lin and C. Liu, Partial regularities of the nonlinear dissipative systems modeling
the flow of liquid crystals, Disc. Contin. Dyn. Syst., A 2 (1996),
1--23.

\bibitem{LW} F. Lin and C. Wang, On the uniqueness of heat flow of harmonic maps and hydrodynamic
flow of nematic liquid crystals, Chinese Annal. Math., 31 (2010),
921--938.

\bibitem{LD} J. Lin and S. Ding, On the well-posedness for the heat flow of harmonic maps and hydrodynamic
flow of nematic liquid crystals in critical spaces, Math. Meth.
Appl. Sciences, DOI: 10.1002/mma.1548.

\bibitem{LNW} C. Liu and N. J. Wakington, Approximation of Liquid
Crystal Flows, SIAM J. Numer. Anal., 37 (2000), 725--741.


%\bibitem{HM} H. Miura, Remark on uniqueness of mild solutions to the Navier-Stokes equations, J. Funct.
%Anal., 218 (2005), 110--129.

\bibitem{JS} J. Serrin, On the regularity of weak solutions of the
Navier-Stokes equations, Arch. Rational Mech. Anal., 9 (1962),
187--195.

%\bibitem{S} E. M. Stein, \textit{Singular integrals and differentiability properties of
%functions}, Princeton, NJ: Priceton University Press, 1971.

\bibitem{SL} H. Sun and C. Liu, On energetic variational approaches in modeling the nematic liquid crystal
flows, Disc. Contin. Dyn. Syst., A 23 (2009), 455--475.

%\bibitem{T} R. Temam, \textit{Navier--Stokes Equations}, North
%Holland, Amsterdam, 1977.

%\bibitem{TH} H. Triebel, \textit{Theory of Function Spaces.}
%Monograph in Mathematics, Vol. 78. Basel: Birkhauser Verglag, 1983.

\bibitem{W08} C. Wang, Heat flow of harmonic maps whose gradients
belong to $L^{n}_{x}L^{\infty}_{t}$, Arch. Rational Mech. Anal., 188
(2008), 309--349.


\bibitem{W} C. Wang, Well-posedness for the heat flow of harmonic
maps and the liquid crystal flow with rough initial data, Arch.
Rational Mech. Anal., 200 (2011), 1--19.

\bibitem{WD} H. Wen and S. Ding, Solutions of incompressible
hydrodynamic flow of liquid crystals, Nonlinear Anal. Real Word
Appl. 12 (2011), 1510--1531.



\end{thebibliography}
\end{document}